\documentclass[10pt,onecolumn,draftcls]{IEEEtran}       %ARXIV
\IEEEoverridecommandlockouts
\usepackage{cite}
\usepackage{amsmath,amssymb,amsfonts,amsthm}
\usepackage[ruled,norelsize]{algorithm2e}
\usepackage{graphicx}
\usepackage{textcomp}
\usepackage{xcolor}
\usepackage{enumerate}
\usepackage{mathtools}
\usepackage{subcaption}
\usepackage{etoolbox}

\providebool{Report}
\setbool{Report}{true}          % TRUE FOR ARXIV; FALSE FOR JOURNAL

\makeatletter
\newcommand{\removelatexerror}{\let\@latex@error\@gobble}
\makeatother

\newcommand{\R}{\mathbb{R}}
\newcommand{\N}{\mathbb{N}}
\newcommand{\M}{\mathcal{M}}

\newtheoremstyle{italicNoExtraSpace}% ⟨name ⟩ 
{3pt}% ⟨Space above ⟩1 
{3pt}% ⟨Space below ⟩1
{\itshape}% ⟨Body font ⟩
{}% ⟨Indent amount ⟩2
{\itshape}% ⟨Theorem head font⟩
{:}% ⟨Punctuation after theorem head ⟩
{.5em}% ⟨Space after theorem head ⟩3
{}% ⟨Theorem head spec (can be left empty, meaning ‘normal’)⟩
\theoremstyle{italicNoExtraSpace}

\newtheorem{assumption}{Assumption}
\newtheorem{theorem}{Theorem}
\newtheorem{problem}{Problem}
\newtheorem{lemma}{Lemma}
\newtheorem{definition}{Definition}
\newtheorem{remark}{Remark}
\newtheorem{corollary}{Corollary}
\newtheorem{example}{Example}

\DeclareMathOperator*{\argmax}{arg\,max}
\DeclareMathOperator*{\argmin}{arg\,min}

\newcommand{\mcc}{\mu^c_{[c]}(t_c)}
\newcommand{\mhatc}{\hat{\mu}_{\delta,[c]}}
\newcommand{\mhat}{\hat{\mu}_{\delta,[c]}}

\newcommand{\xhat}{\hat{x}_{\delta}}
\newcommand{\xhatt}{\hat{x}_{\delta}(t)}
\newcommand{\gc}{g_{[c]}}
\newcommand{\xctc}{x^c_t}
\newcommand{\munormbound}{B}

\newcommand{\mhatd}{\hat{\mu}_{\delta}}

\usepackage{mathtools}
\mathtoolsset{showonlyrefs=true}

\usepackage{cite}

\usepackage{color}

\newcommand{\blue}[1]{{\textcolor{black}{#1}}}      %ARXIV 

\usepackage{marginnote}

\def\BibTeX{{\rm B\kern-.05em{\sc i\kern-.025em b}\kern-.08em
    T\kern-.1667em\lower.7ex\hbox{E}\kern-.125emX}}
\begin{document}

\title{Totally Asynchronous Primal-Dual Convex Optimization in Blocks\\
\thanks{This work was supported by a task order contract from the Air Force Research Laboratory through
Eglin AFB, by ONR under grants N00014-21-1-2495 and N00014-19-1-2543, and by AFOSR under grant FA9550-19-1-0169.}}

\author{
Katherine R. Hendrickson and Matthew T. Hale$^*$\thanks{${}^*$The authors are with the Department of Mechanical and Aerospace Engineering,
Herbert Wertheim College of Engineering, University of Florida, Gainesville, FL~$32611$. Emails:
\texttt{\{kat.hendrickson,matthewhale\}@ufl.edu}. 
}
}

\maketitle

\begin{abstract}
We present a parallelized primal-dual algorithm for solving constrained convex optimization problems. The algorithm is ``block-based,'' in that vectors of primal and dual variables are partitioned into blocks, each of which is updated only by a single processor. We consider four possible forms of asynchrony: in updates to primal variables, updates to dual variables, communications of primal variables, and communications of dual variables. We show that any amount of asynchrony in the communications of dual variables
can preclude convergence, though the other forms of asynchrony are permitted.  
A first-order primal-dual update law is then developed and shown to be robust to these other forms of asynchrony. We next derive convergence rates to a Lagrangian saddle point in terms of the operations agents execute, without specifying any timing or pattern with which they must 
be executed. These convergence rates include an ``asynchrony penalty'' that we quantify and present ways to mitigate. Numerical results illustrate these developments. 
\end{abstract}

\section{Introduction}
A wide variety of machine learning problems can be formalized
as convex programs~\cite{sra12,shalev12,bertsekas15,boyd04}. Large-scale machine learning then requires
solutions to large-scale convex programs, which
can be accelerated through parallelized solvers running
on networks of processors. 
In large networks, it can be difficult to synchronize 
their computations, which
generate new information, and communications, which
share this new information with other processors. Accordingly,
we are interested in asynchrony-tolerant large-scale optimization.

The challenge of asynchrony is that
it causes disagreements
among processors that result from generating and receiving different
information at different times. One way to reduce
disagreements
is through repeated averaging of processors' iterates.
This approach dates back several decades~\cite{tsitsiklis86}, and
approaches of this class include~\cite{nedic09,nedic10,nedic15,duchi11,tsianos12,tsianos12b,zhu11,jaggi14}. 
However, these averaging-based methods require bounded
delays in some form, often through requiring connectedness
of agents' 
communication 
graphs over intervals of a prescribed length~\cite[Chapter 7]{bertsekas89}.
In some applications, delays are outside agents' control, e.g.,
in a contested environment where communications are jammed,
and delay bounds cannot be easily enforced.
Moreover, graph connectivity cannot be easily checked
individual agents, meaning even satisfaction
or violation of connectivity bounds is not readily ascertained. 
In addition, these methods require multiple processors to 
update each decision variable, which can be prohibitive, e.g., in learning 
problems with billions of data points. 

Therefore, in this paper we develop a totally asynchronous parallelized primal-dual
method for solving large constrained convex optimization problems. The term ``totally asynchronous'' dates back to~\cite{bertsekas89} and describes scenarios in which both computations 
and communications are executed without any assumptions on delay bounds. 
By ``parallelized,'' we mean that each decision variable
is updated only by a single processor. As problems grow, this
has the advantage of keeping each processor's computational
burden approximately constant. The decision variables assigned
to each processor are referred to as a ``block,'' and asynchronous
block-based algorithms date back several decades as well~\cite{tsitsiklis86,Bertsekas1983,bertsekas89}.
Those early works solve unconstrained or set-constrained
problems, in addition to select
problems with functional constraints. 
\blue{Recent asynchronous 
block-based 
algorithms have also been developed for some specific 
classes of problems
with set or functional constraints~\cite{cannelli19, cannelli21, peng18, peng15, liu15}}.  
To bring parallelization to arbitrary constrained
problems, we develop a primal-dual approach that does not
require constraints to have a specific form.

Block-based methods have previously been shown to tolerate
arbitrarily long delays in both communications and
computations in some unconstrained problems~\cite{Bertsekas1983,hochhaus18,ubl19},
eliminating the need to enforce
and verify delay boundedness assumptions. 
For constrained
problems of a general form, block-based methods have been
paired with primal-dual algorithms with centralized
dual updates~\cite{hale14,hale17} and/or synchronous primal updates~\cite{koshal2011multiuser}.
To the best of our knowledge, 
arbitrarily asynchronous block-based updates have not been
developed for convex programs of a general form. 
A counterexample in~\cite{hale17} showed that arbitrarily asynchronous
communications of dual variables can preclude convergence, though
that example leaves open the extent to which more limited dual asynchrony
is compatible with convergence.

In this paper, we present 
a primal-dual optimization algorithm that permits
arbitrary asynchrony in primal variables, while
accommodating dual asynchrony to the extent possible.
Four types of asynchrony are possible: (i) asynchrony
in primal computations, 
(ii) asynchrony in communicating primal variables,
(iii) asynchrony in dual computations, and
(iv) asynchrony in communicating dual variables. 
The first contribution of this paper is to show that 
item (iv) is fundamentally problematic.  
Specifically, we show 
that arbitrarily small disagreements among dual variables
can cause primal computations to disagree by arbitrarily large amounts. 
For this reason, we rule out asynchrony in communicating
dual variables. However, we permit all other forms
of asynchrony, and, relative to existing work, this
is the first to permit arbitrarily asynchronous
computations of dual variables in blocks.

The second contribution of this paper is to establish
convergence rates. 
These rates are shown to depend upon problem parameters,
which lets us calibrate their values 
to improve convergence. Moreover, we show that
convergence can be inexact due to dual asynchrony, and thus
the scalability of parallelization comes at the
expense of a potentially inexact solution. We term
this inexactness the ``asynchrony penalty,'' and we give
an explicit bound on it, as well as methods to mitigate it.
Simulation results show convergence of this algorithm 
and illustrate that the asynchrony penalty is mild.

This paper is an extension of the conference paper~\cite{hendrickson20}. 
\blue{This paper extends all previous results on scalar blocks to non-scalar blocks, 
provides bounds on regularization error, provides techniques
to mitigate the asynchrony penalty, and gives a simplified convergence
analysis.}
%and prove that the asynchrony penalty may be completely eliminated for some types of problems. 
%Finally, this paper \blue{provides proofs of results not previously included in the conference paper.}

The rest of the paper is organized as follows.
Section~\ref{sec:background} provides background
and a formal problem statement. Section~\ref{sec:algorithm} presents our asynchronous algorithm. 
Convergence rates are developed in Section~\ref{sec:overall}. 
Section~\ref{sec:numerical} presents simulation results, 
and Section~\ref{sec:concl} concludes.

\section{Background and Problem Statement} \label{sec:background}
\blue{Real-world applications of multi-agent optimization may 
face challenges that 
prevent agents from computing or communicating at specified times or intervals. 
For example, very large networks of processors may face difficulty in synchronizing
all of their clocks, and networks of autonomous agents in a contested environment may face
jammed communications that make information sharing sporadic. 
This asynchrony in computations and communications motivates the development of algorithms that tolerate as much asynchrony as possible. Thus}, we study the following form of optimization problem.

\begin{problem} \label{prob:first}
Given~$f : \R^n \to \R$, ~$g : \R^n \to \R^m$, and~$X \subset \R^n$, asynchronously solve
\begin{align}
\textnormal{minimize }        &f(x) \\
\textnormal{subject to } &g(x) \leq 0 \\
                            &x \in X. \tag*{$\lozenge$}
\end{align} 
\end{problem}

We assume the following about~$f$.

\begin{assumption} \label{as:f}
The objective function~$f$ is twice continuously
differentiable and convex. \hfill $\triangle$
\end{assumption}

We make a similar assumption about
the constraints~$g$.

\begin{assumption} \label{as:g}
For all~$j \in \{1, \ldots, m\}$, the function~$g_j$ is twice continuously differentiable
and convex. 
And~$g$ satisfies Slater's condition, i.e.,
there exists~$\bar{x} \in X$ such
that~$g\big(\bar{x}\big) < 0$. \hfill $\triangle$
\end{assumption}

Assumptions~\ref{as:f} and~\ref{as:g} permit
a wide range of functions to be used, such as all
convex polynomials of all orders. 
We impose the following assumption on
the constraint set. 

\begin{assumption} \label{as:X}
The set~$X$ is non-empty, compact, and convex. 
It can be decomposed via~$X = X_1 \times \cdots \times X_{N_p}$ \blue{where~$N_p$ is the number of agents optimizing over~$x$}. 
\hfill $\triangle$
\end{assumption}

Assumption~\ref{as:X} permits many sets to be used,
such as box constraints, which often arise \blue{in}
multi-agent optimization~\cite{notarnicola2016asynchronous}. 
\blue{This assumption allows agents to project their blocks of the decision variable onto the corresponding part of the constraint set, i.e., for all~$i$, agent~$i$
is able to project its updates onto~$X_i$, which ensures that~$x \in X$ overall. 
This property 
enables a distributed projected update law in which each agent ensures set constraint satisfaction of its block of the decision variable. 
This form of decomposition has been used in~\cite{cannelli19, cannelli21, peng18, liu15, koshal2011multiuser, zhu11} (and other works) for the same purpose.}

We will solve Problem~\ref{prob:first}  
using a primal-dual approach. This allows the problem to be parallelized across many processors by re-encoding constraints through Karush-Kuhn-Tucker (KKT) multipliers. In particular, because the constraints~$g$ couple
the \blue{agents}' computations, they can be difficult to enforce in a distributed way. 
By introducing KKT multipliers to encode
constraints, we can solve an equivalent, higher-dimensional
unconstrained problem. 

An ordinary primal-dual approach would 
find a saddle point of the Lagrangian associated
with Problem~\ref{prob:first},
defined as~$L(x, \mu) = f(x) + \mu^Tg(x)$, where~$\mu \geq 0$. 
That is, one would solve~$\min_{x \in X} \max_{\mu \geq 0} L(x, \mu)$,
and, under Assumptions~\ref{as:f}-\ref{as:X}, this would furnish
a solution to Problem~\ref{prob:first}. 
However,~$L$ is affine in~$\mu$, which implies
that~$L(x, \cdot)$ is concave but not strongly
concave. Strong convexity has
been shown to provide robustness to asynchrony
in minimization problems~\cite{bertsekas89}, and
thus we wish to endow the maximization 
over~$\mu$ with strong concavity. 
We use a Tikhonov regularization~\cite{facchinei2007finite} in~$\mu$
to form
\begin{equation} \label{regL}
L_\delta (x, \mu) = f(x) + \mu^Tg(x) - \frac{\delta}{2} \| \mu \| ^2,
\end{equation}
where~$\delta > 0$ and~$\|\cdot\|$ denotes the Euclidean norm. 
This ensures~$\delta$-strong concavity in~$\mu$. 
\blue{Thus, we will find a saddle point~$(\hat{x}_\delta, \hat{\mu}_\delta)$ of~$L_\delta$,
which is approximately equal to a saddle point of~$L$ and thus approximately
solves Problem~\ref{prob:first}. We bound the error 
introduced by regularizing in Theorem~\ref{lem:regerror} below. }

\blue{One challenge in designing and analyzing an algorithm 
is that~$\mhatd$ is contained
in the unbounded domain~$\R^m_{+}$, which
is the non-negative orthant of~$\R^m$. 
Because this domain is unbounded, \blue{gradients with respect to the dual variable are unbounded}. 
Specifically, dual iterates may not be within a bounded distance of the optimum and hence they may produce
gradients that are arbitrarily large. 
%which makes convergence analysis challenging because dual iterates may not be within a bounded distance of the optimum
%and hence they may produce
%gradients that are arbitrarily large. 
To remedy this problem, we will confine dual variables to a set~$\mathcal{M} \subset \mathbb{R}^m_{+}$, defined as
%convex set~$\mathcal{M}$:
\begin{equation} \label{eq:Mdef}
    \mathcal{M} := \big\lbrace \mu \in \R^m_+ : \| \mu \|_1 \leq \munormbound \big\rbrace, \quad \munormbound := \frac{f(\bar{x}) - f^*}{\min\limits_{1 \leq j \leq m}  -g_j(\bar{x})},
\end{equation}
where~$\bar{x}$ is any Slater point.}
\blue{Here,~$f^*$ denotes the optimal objective
function value over~$X$ (but without~$g$), though any lower-bound for this value will suffice. For example, if~$f$ is non-negative, then
one can substitute~$0$ in place of~$f^*$.
We will show below in Lemma~\ref{lem:mubound}  that using~$\mathcal{M}$ does not affect
the final answer that is computed. 
}

Instead of regularizing with respect to the primal variable~$x$, we impose the following assumption
in terms of the Hessian\blue{~$H(x,\mu) := \nabla_{x}^2 L_{\delta}(x, \mu)$.} When convenient, we suppress
the arguments~$x$ and~$\mu$ and simply write~$H$.

\begin{assumption}[Diagonal Dominance] \label{as:diagonal}
The Hessian matrix~$H=\nabla^{2}_{x}L_{\delta}(x,\mu)$ is~$\beta$-diagonally dominant \blue{for all~$\mu \in \M$}. That is,~$|H_{ii}|-\beta \geq \sum_{\substack{ j=1 \\ j \neq i}}^n |H_{ij}|, \forall i = 1, \ldots, n.$ \hfill $\triangle$
%\begin{equation}
%|H_{ii}|-\beta \geq \sum_{\substack{ j=1 \\ j \neq i}}^n |H_{ij}|
%\qquad \textnormal{for all } i = 1, \ldots, n. 
%\tag*{$\triangle$}
%\end{equation}
\end{assumption}

\blue{
It has been observed in the literature that Assumption~\ref{as:diagonal} or a similar variant
of diagonal dominance is necessary to ensure the convergence of totally asynchronous
algorithms~\cite[Section 6.3.2]{bertsekas89}. 
}
If this assumption does not hold, the Lagrangian can be regularized with respect to~$x$ to 
help provide~$H$'s diagonal dominance. 
\blue{If this is done, there may be cases in which the 
regularization parameter required to satisfy Assumption~\ref{as:diagonal} is large, and
this can introduce large regularization errors, which can be undesirable; 
see~\cite{koshal2011multiuser} for bounds on regularization error when both primal
and dual regularizations are used.}
Fortunately, numerous problems 
satisfy
this assumption without regularizing in~$x$~\cite{greene06}, and,
for such problems, we proceed without regularizing in~$x$
to avoid unnecessarily introducing regularization error. 
\blue{Diagonal dominance 
has been shown to arise
in sum of squares problems~\cite{ahmadi17}, linear systems with sparse graphs~\cite{cai19}, matrix scaling and balancing~\cite{cohen17}, and quadratic programs~\cite{ubl19}.}
%When applying our asynchronous algorithm, however, diagonal dominance improves convergence. This follows examples in the
%literature that illustrate the benefit of diagonal dominance in optimization
We show in Section~\ref{sec:numerical} that diagonal dominance improves
convergence of our algorithm, and this is in line with existing 
algorithms~\cite{cohen17, ahmadi17,ahmadi19,frommer91,cai19,zhang19}.

\blue{Using the definition of~$\mathcal{M}$ in~\eqref{eq:Mdef} and Assumption~\ref{as:diagonal}, we now observe that~$\mathcal{M}$ contains the optimum~$\mhatd$. \begin{lemma} \label{lem:mubound}
Let Assumptions~\ref{as:f}-\ref{as:diagonal} hold. % and
%let~$\bar{x}$ be a Slater point of~$g$. Set~$f^* := \min\limits_{x \in X} f(x)$
%and~$\munormbound := \frac{f(\bar{x}) - f^*}{\min\limits_{1 \leq j \leq m}  -g_j(\bar{x})} $. 
Then
%\begin{equation*}
$\mhatd \in \mathcal{M}$.% := \big\lbrace \mu \in \R^m_+ : \| \mu \|_1 \leq \munormbound \big\rbrace$.
%\end{equation*}
\end{lemma}
\emph{Proof:} 
Follows Section~II-C in~\cite{hale15}. \hfill $\blacksquare$}

\blue{We now present the following saddle point problem that will be the focus of the rest of the paper.}

\begin{problem} \label{prob:second}
Let Assumptions~\ref{as:f}-\ref{as:diagonal} hold and fix~$\delta > 0$. 
For~$L_{\delta}$ defined in~\eqref{regL}, asynchronously compute
\begin{equation} 
\big(\xhat, \mhatd ) := \argmin_{x \in X} \argmax_{\mu \in \mathcal{M}}  L_\delta (x, \mu). \tag*{$\lozenge$}
\end{equation}
\end{problem}

The strong convexity of~$L_{\delta}(\cdot, \mu)$ and strong concavity
of~$L_{\delta}(x, \cdot)$ imply that~$(\xhat, \mhatd)$ is unique. 
\blue{
Due to regularizing, the solution to Problem~\ref{prob:second} may not 
equal that of Problem~\ref{prob:first}, and regularization could also introduce constraint violations. We next bound both regularization error in solutions
and constraint violations in terms of the regularization
parameter~$\delta$.}
\blue{\begin{theorem}\label{lem:regerror}
Let Assumptions~\ref{as:f}-\ref{as:diagonal} hold. 
Let~$(\hat{x}, \hat{\mu})$ denote a saddle point of~$L$ (without regularization applied). 
Then the regularization error introduced by the Tikhonov regularization
in~\eqref{regL} is bounded by~$\|\hat{x}_\delta - \hat{x}\|^2 \leq \frac{\delta}{\beta}\munormbound^2$, 
where~$\delta$ is the regularization parameter and~$\munormbound$ is defined in Lemma~\ref{lem:mubound}. Furthermore, possible constraint violations 
are bounded via~$g_j(\hat{x}_\delta)\leq M_j\munormbound\sqrt{\frac{\delta}{\beta}}$,
where~$M_j := \max_{x \in X} \|\nabla g_j(x)\|$.
\end{theorem}
\emph{Proof:} 
See Appendix~\ref{app:regerror}. \hfill $\blacksquare$}

\blue{
The error in solutions is~$O(\delta)$ and the potential constraint
violation is~$O(\sqrt{\delta})$, and thus both can be made arbitrarily
small. Moreover, if it is essential that a feasible point be computed,
then, for all~$j$, one can replace the constraint~$g_j(x) \leq 0$
with~$\tilde{g}_j(x) = g_j(x) - M_j\munormbound\sqrt{\frac{\delta}{\beta}} \leq 0$,
which will ensure the generation of a feasible point. 
And Slater's condition in Assumption~\ref{as:g} implies that, for
sufficiently small~$\delta$, there exist points
that satisfy~$\tilde{g}_j$. 
}

\section{Asynchronous Primal-Dual Algorithm} \label{sec:algorithm}
Solving Problem~\ref{prob:second} asynchronously requires an update law that 
we expect to be robust to asynchrony and simple to implement in
a distributed way. In this context, first-order
gradient-based methods offer some degree of inherent
robustness, as well as computations that are simpler than
other methods, such as Newton-type methods. 
We apply a projected gradient method to both the primal and dual variables, based on the 
seminal Uzawa algorithm~\cite{arrow58}. \blue{Recall
that, given some~$x(0)$ and~$\mu(0)$, at iteration~$k+1$ 
the Uzawa algorithm computes
the primal update,~$x(k+1)$, and dual update,~$\mu(k+1)$, using
\begin{align}
x(k+1) &= \Pi_X [x(k) - \gamma \nabla_x L_\delta \big(x(k),\mu(k)\big)] \label{eq:uzawax}\\ 
\mu(k+1) &= \Pi_\mathcal{M} [\mu(k) + \rho \nabla_\mu L_\delta \big(x(k),\mu(k)\big)], \label{eq:uzawamu}
\end{align}
where~$\gamma, \rho > 0$ are stepsizes,~$\Pi_X$ is the Euclidean projection onto~$X$, and~$\Pi_\mathcal{M}$ is the Euclidean projection onto~$\mathcal{M}$.}

\subsection{Overview of Approach}
\blue{The Uzawa algorithm is centralized, and} we will decentralize~\eqref{eq:uzawax} and~$\eqref{eq:uzawamu}$ among a number of agents while allowing them to 
generate and share information
as asynchronously as possible. We consider~$N$ agents indexed over~$i \in \mathcal{I}:=\lbrace 1, \ldots , N\rbrace$. We also define
the sets~$\mathcal{I}_p:=\lbrace 1, \ldots , N_p \rbrace$ and~$\mathcal{I}_d:=\lbrace 1, \ldots , N_d \rbrace$, 
%where~$N_p$ is the number of ``primal agents'' (agents that update primal variables), and~$N_d$ is the number of ``dual agents'' 
%(agents that update dual variables). 
where~$N_p+N_d = N$. 
The set~$\mathcal{I}_p$ contains indices of ``primal agents'' that update primal blocks (contained in~$x$), while~$\mathcal{I}_d$ contains indices of 
``dual agents'' 
that update dual blocks (contained in~$\mu$)\footnote{
Although the same index may be contained in both~$\mathcal{I}_p$ and~$\mathcal{I}_d$, we define
the sets in this way to avoid non-consecutive numbering of primal agents and dual agents, which would
be cumbersome in the forthcoming analysis. The meaning of an index will always be made unambiguous
by specifying whether it is contained in~$\mathcal{I}_p$ or~$\mathcal{I}_d$. 
}. 
Thus,~$x \in \R^n$ is divided into~$N_p$ blocks and~$\mu \in \R^m$ into~$N_d$ blocks.

Primal agent~$i$ updates the~$i^{th}$ primal block,~$x_{[i]}$, and dual agent~$c$ updates the~$c^{th}$ dual block,~$\mu_{[c]}$. Let~$n_i$ denote the length of primal agent~$i$'s block and~$m_c$ the length of dual agent~$c$'s block. Then~$n= \sum_{i=1}^{N_p}n_i$ and~$m= \sum_{c=1}^{N_d}m_c$. The block of constraints~$g$ 
that correspond to~$\mu_{[c]}$ is denoted by~$g_{[c]} : \R^n \to \R^{m_c}$ and each dual agent projects its 
computations onto a set~$\mathcal{M}_{c}$ derived from Lemma~$\ref{lem:mubound}$, namely
%\begin{equation}
    $\mathcal{M}_{c} = \big\{\nu \in \mathbb{R}_{+}^{m_c} : \| \nu \|_1 \leq \munormbound\big\}$.%\frac{f(\bar{x}) - f^*}{\min\limits_{1 \leq j \leq m} -g_j(\bar{x})}\right\}.
%\end{equation}
%The set~$\mathcal{M}_1 \times \cdots \times \mathcal{M}_{N_d}$ is larger than~$\mathcal{M}$, but 
%suffices for our developments. 

Using a primal-dual approach, there are four behaviors
that could be asynchronous: (i) computations of primal variables, (ii) communications
of the values of primal variables, (iii) computations of
dual variables, and (iv) communications of the 
values of dual variables. \blue{In all cases, we assume that communications arrive in finite time and are received in the order they were sent.} We examine these four behaviors here:

\paragraph*{(i) Computations of Updates to Primal Variables} When parallelizing~\eqref{eq:uzawax} across the~$N_p$ primal agents, 
we index all primal agents' computations
using the same iteration counter,~$k \in \N$. 
However, they may compute and communicate at different times and they do not necessarily do either at all~$k$. 
The subset of times at which primal agent~$i \in \mathcal{I}_p$ computes an update is denoted by~$K^i \subset \N$. 
For distinct~$i, j \in \mathcal{I}_p$, we allow~$K^i \neq K^j$. These sets
are used only for analysis and need not be known to agents. 

\paragraph*{(ii) Communications of Primal Variables} Primal variable communications are also totally asynchronous. A primal block's current value may or may not be sent to other primal and dual agents that need it at each time~$k$. Thus, one agent may have onboard an old value of a primal block
computed by another agent. 
\blue{We use~$\mathcal{N}_i \subset \mathcal{I}_p$ to denote the set of indices of primal agents whose
decision variables are needed for agent~$i$'s computations. 
Formally,~$j \in \mathcal{N}_i$ if and only if~$\nabla_{x_i} L_{\delta}(x, \mu)$
explicitly depends on~$x_{[j]}$. The set~$\mathcal{N}_i$ is referred to 
as agent~$i$'s \emph{essential neighbors}, and only agent~$i$'s essential neighbors need
to communicate to agent~$i$. In particular, primal communications are not all-to-all. 
We use~$\tau^i_j(k)$ to denote 
the time at which primal agent~$j$ originally computed the value of~$x_{[j]}$
stored onboard primal agent~$i$ at time~$k$.
We use~$\sigma^c_j(k)$ to denote the time at which primal agent~$j$
originally computed the value of~$x_{[j]}$ stored onboard 
dual agent~$c$ at time~$k$. 
These functions are used only for analysis, i.e., agents do not need to know the values of~$\tau^i_j(k)$ or~$\sigma^c_j(k)$.}

\blue{We impose the following assumption.
\begin{assumption}[Primal Updates and Communications] \label{as:times}
For all~$i \in \mathcal{I}_p$, the set~$K^i$ is infinite. If~$\{k_n\}_{n \in \N}$ is an increasing set
of times in~$K^i$, then~$\lim_{n \to \infty} \tau^j_i(k_n) = \infty$ for all~$j \in \mathcal{I}_p$ such that~$i \in \mathcal{N}_j$
and~$\lim_{n \to \infty} \sigma^c_i(k_n) = \infty$ for all~$c \in \mathcal{I}_d$ such that~$x_{[i]}$ is constrained by~$g_{[c]}$. \hfill $\triangle$
\end{assumption}
This simply ensures that, for all~$i \in \mathcal{I}_p$, primal agent~$i$ never stops computing or 
communicating, though delays can be arbitrarily large.}

\paragraph*{(iii) Computations of Updates to Dual Variables} 
Dual agents wait for every primal agent's updated block before computing an update to a dual variable. 
Dual agents may perform computations at different times because they may receive updates to primal blocks at different times. In some cases, a dual agent may receive multiple updates from a subset 
of primal agents prior to receiving all required primal updates. In this case, only the most recently received update from a primal agent will be used in the dual agent's computation.
For all~$c \in \mathcal{I}_d$, dual agent~$c$ keeps an iteration count~$t_c$ to track the number of updates it has completed. 
%To simplify notation, the time variable~$t$ used by the dual agents is a vector consisting of~$t_{c}$ for all~$c \in \mathcal{I}_d$.

\paragraph*{(iv) Communications of Updated Dual Variables} 
\blue{Previous work~\cite[Section VI]{hale17} has shown that allowing primal agents to disagree arbitrarily about dual variables
can preclude convergence. In particular, that work provides an example problem in which such disagreements lead to oscillations 
in the primal variables
that do not decay with time, and
thus agents do not even converge.}
This is explained by the following: fix~$\mu^1, \mu^2 \in \mathcal{M}$. Then a primal agent with~$\mu^1$ onboard
is minimizing~$L_{\delta}(\cdot, \mu^1)$, while a primal agent with~$\mu^2$ onboard is minimizing~$L_{\delta}(\cdot, \mu^2)$. 
If~$\mu^1$ and~$\mu^2$ can be arbitrarily far apart, then it is not surprising that
the minima of~$L(\cdot, \mu^1)$ and~$L(\cdot, \mu^2)$ are arbitrarily far apart, which is what is observed in~\cite[Section VI]{hale17}. 
One may then conjecture that small disagreements in dual variables lead to small
distances between these minima. We next show that this is false. 

\begin{theorem} \label{thm:disagree}
Fix any~$\epsilon > 0$ and any~$L > \epsilon$. Then, under
Assumptions~\ref{as:f}-\ref{as:diagonal}, there exists a problem and points~$\mu^1, \mu^2 \in \mathcal{M}$ such
that~$\|\mu^1 - \mu^2\| < \epsilon$ 
and~$\|\hat{x}_1 - \hat{x}_2\| > L$,
where~$\hat{x}_1 = \argmin_{x \in X} L_{\delta}(x, \mu^1)$
and~$\hat{x}_2 = \argmin_{x \in X} L_{\delta}(x, \mu^2)$. 
\end{theorem}
\emph{Proof:} See Theorem~1 in the preliminary version
of this work~\cite{hendrickson20arxiv}. \hfill $\blacksquare$ 

\blue{The counterexample
in~\cite[Section VI]{hale17} shows that primal agents need not even converge
if they use different dual variables that take them to minima that are far apart.
Theorem~\ref{thm:disagree} shows that arbitrarily small disagreements in the values
of dual variables can drive primal agents' computations to points that are arbitrarily far apart. 
Thus, by combining these two results, 
we see that \emph{any} disagreement in the dual variables can preclude convergence.
} 
%Even limited asynchrony can lead to small disagreements in dual
%variables, and, in light of the above discussion, this can cause primal agents' computations to reach points
%that are arbitrarily far apart. 
%
%
\blue{
Therefore, 
%we will develop an algorithm that proceeds
%with all primal agents agreeing on the value of~$\mu$ while still allowing 
%dual
%computations to be divided among dual agents. 
primal agents are allowed to operate totally 
asynchronously, but their computations must use the same dual
variable, formalized as follows.
\begin{assumption}[Dual Communications] \label{as:dualcomm}
Any transmission sent from primal agent~$i$ to primal agent~$j$ while they both have~$\mu(t)$ onboard is only used by primal agent~$j$ in its own computations 
if it is received before the next dual update. \hfill $\triangle$
\end{assumption}
}

\blue{
However, we emphasize that this \emph{does not} mean that the values
of dual variables must be communicated synchronously. 
Instead, primal
agents can use any method to ensure that their computations use
the same dual variables. 
For example, when dual agent~$c$ sends its updated dual block to primal agents, it 
can also send its iteration count~$t_{c}$. Primal agents can use these~$t_c$ values
to annotate which version of~$\mu$ is used in their updates, e.g.,
by appending the value of~$t_c$ for each~$c$ to the end of the vector of primal variables they communicate. 
To ensure that further primal updates rely upon the same dual value, 
other primal agents will
disregard any received primal updates that are annotated with an old iteration count 
for any block of the dual variable. 
%Then, if a primal agent has a dual variable with an iteration counter of~$t_c$ for block~$c$, this agent will ignore a primal
%update it receives if it was computed with an iteration counter of~$t_c\!-\!1$ (or earlier) for dual block~$c$. 
%%Similar to~\cite{hale17}, this process is formalized in Assumption~\ref{as:dualcomm} below.
}

\blue{
The algorithm we present is unchanged if any other method is used to ensure
that the primal agents use the same dual variable in their computations, e.g., 
primal agents may send each other acknowledgments of a new dual variable prior to computing new iterations.
}

\blue{We also note that for many problems, dual agents may not be required to communicate with all primal agents. 
%Our algorithm is presented in the most general manner so that it may be applied to the largest problem set. 
Problems with a constraint function~$g_c$ that depends only on a subset of primal variables will result in a dual entry~$\mu_c$ that needs to be sent only to the
corresponding subset of primal agents, which we illustrate next. 
}

\begin{example}[Dual-to-Primal Communications] \label{ex:onlyone}
\blue{Consider a problem with any objective function~$f$,~$x \in \mathbb{R}^3$, and constraints~$g: \R^3 \to \R^3$ given by
\begin{equation}
    g_1(x) = x_1 + x_2 - b_1, \quad g_2(x) = x_2 - b_2, \quad g_3(x) = x_3 - b_3,
\end{equation}
where~$b_1$,~$b_2,$ and~$b_3$ are some constants.
The regularized Lagrangian associated with this problem is}
\ifbool{Report}{\begin{equation}
L_{\delta}(x, \mu) = f(x) + \mu_1(x_1 + x_2 - b_1)
+ \mu_2(x_2 - b_2)
+ \mu_3(x_3 - b_3) - \frac{\delta}{2} \| \mu \|^2 .
\end{equation}}{\blue{\begin{multline}
L_{\delta}(x, \mu) = f(x) + \mu_1(x_1 + x_2 - b_1)
+ \mu_2(x_2 - b_2) \\
+ \mu_3(x_3 - b_3) - \frac{\delta}{2} \| \mu \|^2 .
\end{multline}}}
\blue{We observe that, among~$\mu_1$,~$\mu_2$, and~$\mu_3$, $\nabla_{x_1}L_{\delta}(x, \mu)$
depends only on~$\mu_1$. Similarly,
$\nabla_{x_2}L_{\delta}(x, \mu)$ depends
only on~$\mu_1$ and~$\mu_2$, 
and~$\nabla_{x_3}L_{\delta}(x, \mu)$ depends
only on~$\mu_3$. 
%In the Lagrangian~$L_{\delta}$ associated with
%this problem, we will have~$\mu \in \mathbb{R}^3$. 
%The dual variable~$\mu_1$ will be used to enforce~$g_1$,~$\mu_2$ 
%will enforce~$g_2$, and~$\mu_3$ will enforce~$g_3$. 
Therefore, only primal agents computing~$x_1$ and~$x_2$ 
would need to receive~$\mu_1$, a primal agent computing~$x_2$ 
would need~$\mu_2$, and only a primal agent computing~$x_3$ would need~$\mu_3$. 
For primal blocks that are scalar values, this leads to the required dual-to-primal communications shown in Figure~\ref{fig:dualcomms1}.}

\blue{Furthermore, the primal and dual variables may instead be divided into the blocks shown in Figure~\ref{fig:dualcomms2}, 
leading to an even simpler communication requirement. In that case, only Primal Agent~1 needs updates from Dual Agent~1 and
only Primal Agent~2 needs updates from Dual Agent~2.}

\begin{figure}[h!]
         \centering
         \includegraphics[width=5cm]{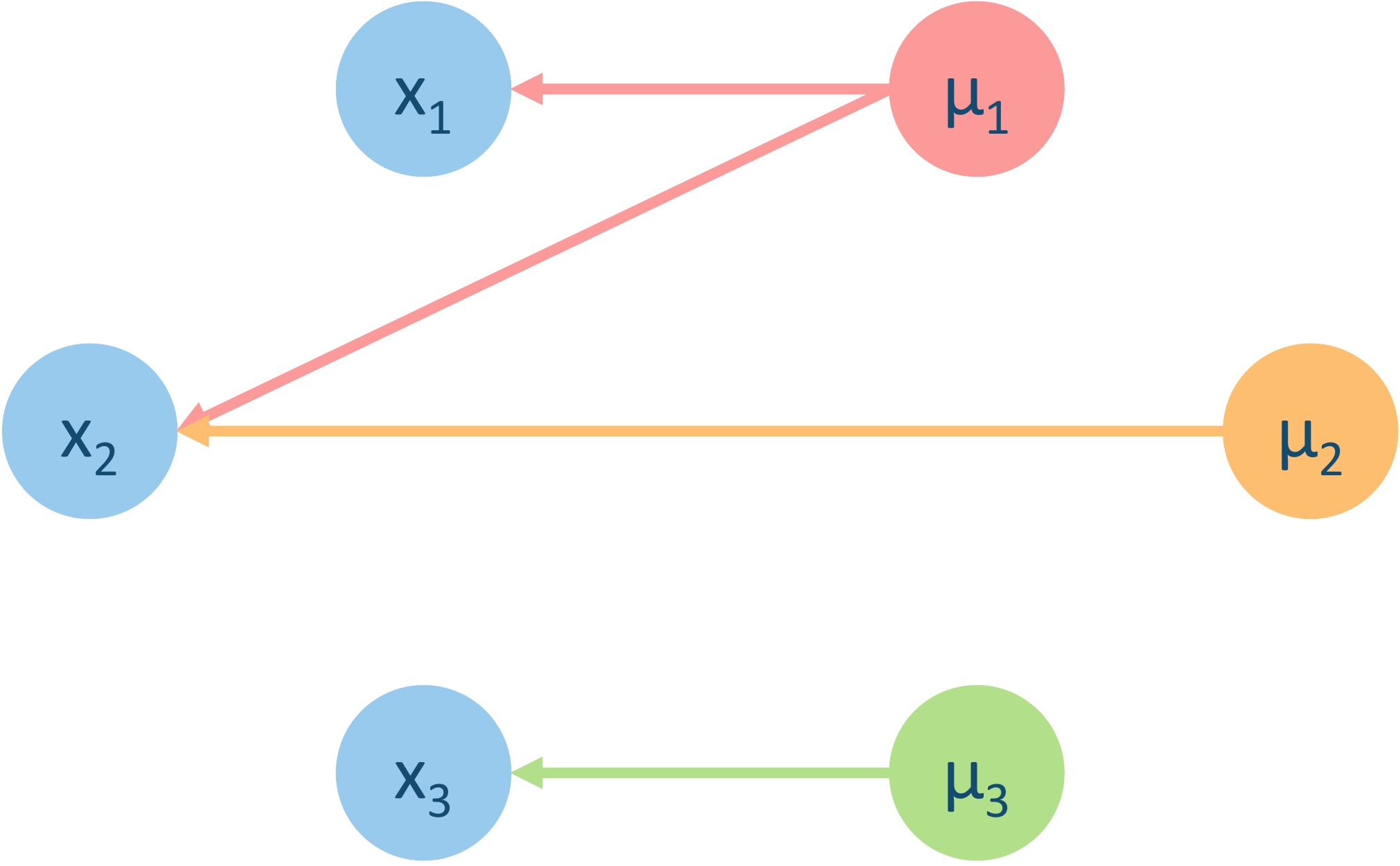}
         \caption{\blue{Required dual-to-primal communications in Example~\ref{ex:onlyone} with scalar blocks. 
         This illustrates that some constraint formulations will only require dual agents to communicate to a subset of primal agents.}}
         \label{fig:dualcomms1}
\end{figure}

\begin{figure}[h!]
         \centering
         \includegraphics[width=5cm]{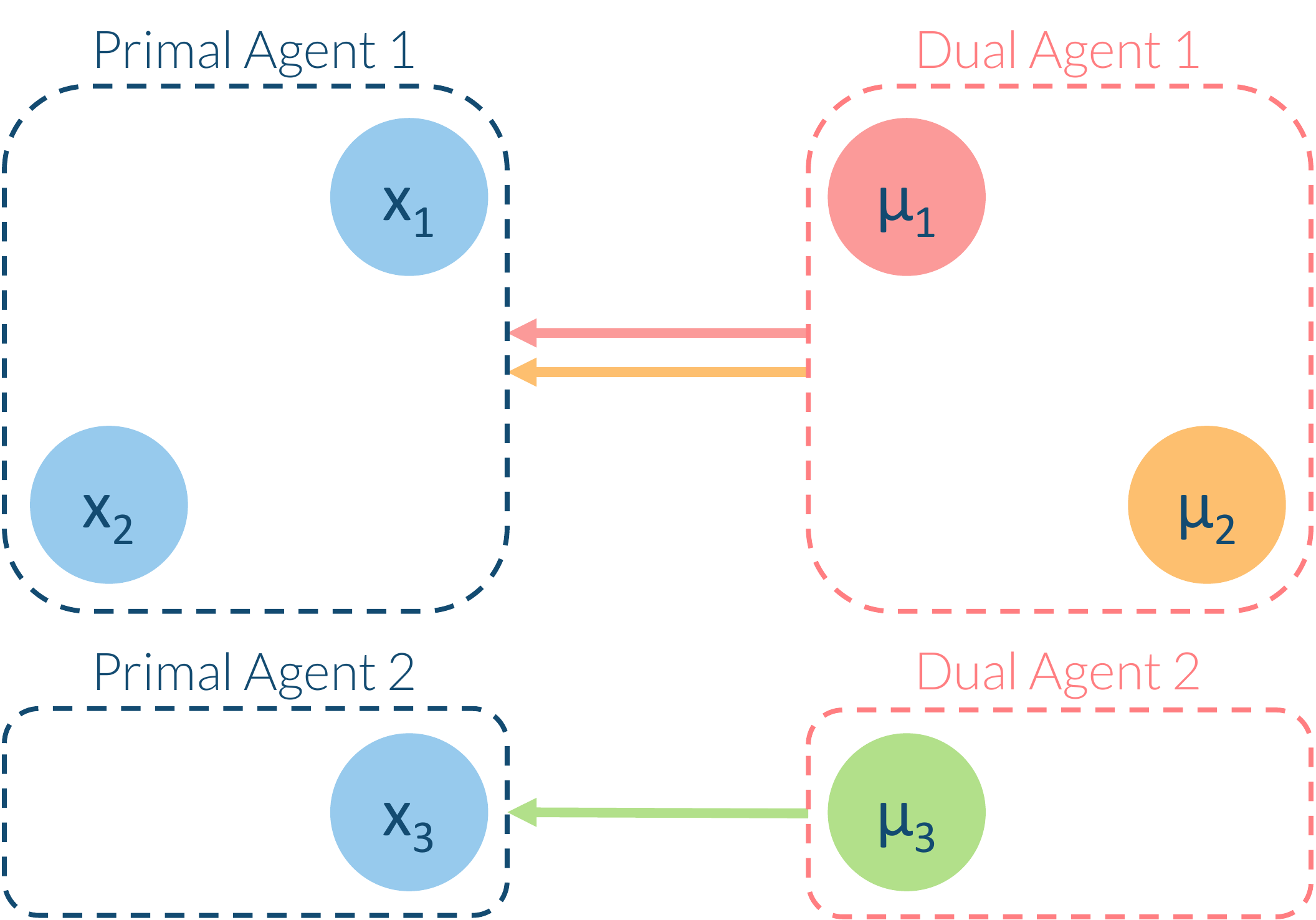}
         \caption{\blue{Required dual-to-primal communications in Example~\ref{ex:onlyone} when separating the primal and dual variables into non-scalar blocks. By dividing blocks according to constraints, required dual-primal communications may be reduced even further.}}
         \label{fig:dualcomms2}
\end{figure} 
\end{example}

\blue{
The principle illustrated by this example
is that, by using non-scalar blocks and exploiting the structure of a 
problem, the required dual
communications may be significantly reduced.
Specifically, a dual block~$\mu_{[c]}$ only needs to be sent to
the primal agents whose decision variables appear
in the block of constraints~$g_{[c]}$. This is reflected
in Step~7 our statement of Algorithm~\ref{alg:2} below. 
}

\subsection{Glossary of Notation in Algorithm~\ref{alg:2}}
The following glossary contains the notation used in our algorithm statement:
\begin{itemize}
%\item[$g_{[c]}(x)$] The~$c^{th}$ block of the constraint function,~$g$, evaluated at~$x$.
\item[$k$] The iteration count used by all primal agents.
\item[$K^{i}$] The set of times at which primal agent~$i$ computes updates.
\item[$\nabla_{x_{[i]}}$] The derivative with respect to the~$i$-th block of~$x$. That is,~$\nabla_{x_{[i]}}:= \frac{\partial}{\partial x_{[i]}}$.
\item[$\mathcal{N}_i$] Essential neighborhood of primal agent~$i$. 
\item[${\mathcal{I}_d}$] Set containing the indices of all dual agents.
\item[${\mathcal{I}_p}$] Set containing the indices of all primal agents.
\item[$\sigma^{c}_{j}(k)$] \blue{Time at which primal agent~$j$ originally computed the value of~$x_{[j]}$ onboard 
dual agent~$c$ at time~$k$.} 
\item[$\tau^{i}_{j}(k)$] \blue{Time at which primal agent~$j$ originally computed the value of~$x_{[j]}$ onboard 
primal agent~$i$ at time~$k$.} Note that~$\tau^{i}_{i}(k)=k$ for all~$i \in \mathcal{I}_p$.
\item[$t$] The vector of dual agent iteration counts. The~$c^{th}$ entry,~$t_c$, is the iteration count for dual agent~$c$'s updates.
%\item[$t_{c}$] The iteration count for dual agent~$c$'s updates. This is sent along with~$\mu^{c}_{c}$ to all primal agents.
\item[$x^{i}_{[j]}$] Primal or dual agent~$i$'s value for the primal block~$j$, which is updated/sent by primal agent~$j$. If agent~$i$ is primal, it is indexed by both~$k$ and~$t$; if agent~$i$ is dual, it is indexed only by~$t$.
%\item[$x^{i}_{[i]} (k;t)$] Primal agent~$i$'s value for its own primal block~$i$ at primal time~$k$, calculated with the value of~$\mu$ from dual update~$t$. \kh{i also think this item might be cut from the glossary, esp. since there is no corresponding dual one.}
%\item[$x^*(t)$] The fixed point of~$h(x)=\Pi_{X}\left[x - \gamma\nabla_{x}L_{\delta}(x,\mu(t))\right]$ with respect to a fixed~$\mu(t)$.
%\item[$x^c_t$] Abbreviation for $x^c(t)$, which is dual agent~$c$'s copy of the primal vector at time~$t$.
\item[$\xhat$] The primal component of the saddle point of~$L_{\delta}$. Part of the optimal solution pair~$(\xhat , \mhatd)$.
%\item[$\xhatt$] Given~$\mu(t)$,~$\xhatt = \argmin_{x \in X} L_{\delta}(x, \mu (t))$.
\item[$\mu^{c}_{[d]}$] Primal or dual agent~$c$'s copy of dual block~$d$, which is updated/sent by dual agent~$d$.
\item[$\mhatd$] The dual component of the saddle point of~$L_{\delta}$,~$(\xhat , \mhatd)$.
%\item[$\mhatc$] The~$c^{th}$ block of~$\mhatd$.
%\item[$\mathcal{M}_{c}$] The set~$\{\nu \in \mathbb{R}_{+}^{m_c} : \| \nu \|_1 \leq\frac{f(\bar{x}) - f^*}{\min_{j} -g_j(\bar{x})}\}$. This uses the upper bound in Lemma~\ref{lem:mubound} to project individual blocks of~$\mu$. 
\end{itemize}

\subsection{Statement of Algorithm}
We now state the asynchronous primal-dual algorithm. 

\begin{figure}[!h]
 \removelatexerror
  \begin{algorithm}[H]
\caption{}
Step 0: Initialize all primal and dual agents with~$x(0)\in X$ and~$\mu(0)\in \mathcal{M}$. 
Set~$t=0 \in \R^{N_d}$ and~$k=0 \in \N$. \\
Step 1: \blue{For all~$i\in\mathcal{I}_p$ and all~$c\in\mathcal{I}_d$, if primal agent~$i$ receives a dual variable update from dual agent~$c$, it sets}
\begin{equation*}
\mu^{i}_{[c]}(t_c) = \mu^{c}_{[c]}(t_{c}).
\end{equation*}

Step 2: For all~$i\!\in\!\mathcal{I}_p$, if~$k \!\in\! K^{i}$, primal agent~$i$ executes
\begin{align*}
\!\!x^{i}_{[i]}(k\!+\!1;t) \!&=\! \Pi_{X_{i}} [x^{i}_{[i]}(k;t) \!-\! \gamma\nabla_{x_{[i]}}L_{\delta}(x^{i}(k;t),\mu^i(t))].
\end{align*}
If~$k \!\notin\! K^{i}$, then $x^{i}_{[i]}(k\!+\!1;t) \!=\!x^{i}_{[i]}(k;t)$.

Step 3: For all~$i\in\mathcal{I}_p$ and all~$j\in \mathcal{N}_{i}$,    
\begin{align*}\!\!x^{i}_{[j]}(k\!+\!1;t)\!=\! &\begin{cases} 
                              x^{j}_{[j]}(\tau^{i}_{j}(k\!+\!1);t) & \!\!\!i$ receives~$x^j_{[j]}\!$ at~$k\!+\!1 \\
                              x^{i}_{[j]}(k;t) & \!\!\!$otherwise$
                              \end{cases}                    
\end{align*}

\blue{Step 4: For all~$i\!\in\!\mathcal{I}_p$, primal agent~$i$ may send~$x^{i}_{[i]}(k+1;t)$ to any primal or dual agent. Due to communication delays, it may not be received for some time.} Set~$k:=k+1$.\\
Step 5: For~$c \in \mathcal{I}_d$ and~$i\in\mathcal{I}_p$, if dual agent~$c$ receives an update from primal agent~$i$ computed with dual update~$t$, it sets
\begin{align*}
x^{c}_{[i]}(t_{c})\!&=x^{i}_{[i]}(\blue{\sigma^{c}_{i}}(k);t).
\end{align*}
Otherwise,~$x^{c}_{[i]}(t_{c})$ remains constant. \\
Step 6: For~$c \in \mathcal{I}_d$, if dual agent~$c$ has received an update from every primal agent constrained by~$g_{[c]}$ 
%\kh{do they really need updates from all primal agents, or should we limit it here to all primal agents that are constrained by~$g_{[c]}$?}
that was computed with the latest dual iteration~$t$, dual agent~$c$ executes
\begin{equation*}
\mu^{c}_{[c]}(t_{c}+1) = \Pi_{\mathcal{M}_{c}}[\mu^{c}_{[c]}(t_{c}) + \rho\frac{\partial \blue{L_\delta}}{\partial \mu_{[c]}}(x^{c}(t_{c}),\mu^{c}(t_{c}))]. 
\end{equation*}

Step 7: If dual agent~$c$ updated in Step 6, it sends~$\mu^{c}_{[c]}(t_c\!+\!1)$ to all primal agents \blue{that are constrained by~$g_{[c]}$. Due to asynchrony, it may not be received for some time.} 
Set~$t_{c}:=t_{c}\!+\!1$.\\
Step 8: Return to Step 1.
\label{alg:2}
\end{algorithm}
\end{figure}

\ifbool{Report}{\newpage}{}
\section{Overall Convergence and Reducing the Asynchrony Penalty}\label{sec:overall}
\blue{In this section, we present our main convergence result and strategies for reducing the asynchrony penalty, which is an error term in that
result that is due to asynchronous operations.}
\blue{
First, 
let~$H(x, \mu) =\nabla^{2}_{x}L_{\delta}(x,\mu)$
and choose the primal stepsize~$\gamma > 0$ to satisfy
\begin{equation} \label{eq:gamma}
\gamma < \frac{1}{\max\limits_i\max\limits_{x \in X} \max\limits_{\mu \in \mathcal{M}} \sum_{j=1}^n |H_{ij}(x, \mu)|} .
\end{equation}
%\end{definition}
Recall that~$x$ and~$\mu$ both take values in compact sets (cf. Assumption~\ref{as:X} and Lemma~\ref{lem:mubound}), and thus 
%These compact sets are particularly useful as they bound both~$x$ and~$\mu$, leading to the boundedness of various functions. It is useful in the definition above, as it implies that 
each entry~$H_{ij}$ is bounded. In particular, 
the upper bound on~$\gamma$ is positive.}

%The primal-only convergence rate established in Appendix~\ref{app:primal} can be computed by
%leveraging results in~\cite{hale17} in terms of the number
%of operations the primal agents have completed (counted in the appropriate
%sequence). 
The main convergence result for Algorithm~\ref{alg:2} is in terms of the number
of operations that agents have executed, counted in a specific order as follows. 
%We count operations as follows. 
Upon the first primal agent's receipt of a
dual variable with iteration vector~$t$, we set~$\textnormal{ops}(k,t) = 0$. Then, after
all primal agents have computed an update to their decision variable with~$\mu(t)$ and
sent it to and had it received by all other primal agents
in their essential neighborhoods, say by time~$k'$, we increment~$\textnormal{ops}$
to~$\textnormal{ops}(k',t) = 1$. After~$\textnormal{ops}(k',t) = 1$,
we then wait until all primal agents have subsequently computed
a new update (still using the same dual variable indexed with~$t$) and it has been 
sent to and received by all primal agents' essential neighbors. 
If this occurs at time~$k''$, then
we set~$\textnormal{ops}(k'',t) = 2$, and then this process
continues. If at some time~$k'''$, primal agents receive an updated~$\mu$ (whether just a single dual agent sent an update or multiple dual agents send updates) with an iteration vector of~$t'$, then the count would begin again with~$\textnormal{ops}(k''',t') = 0$.

\subsection{\blue{Main Result}}
\ifbool{Report}{We now present our main result on the convergence of~$x^i(k; t)$ to~$\xhat$. Recall that~$\delta$ is the dual regularization parameter,~$\gamma$ is the primal 
stepsize given in~\eqref{eq:gamma}, and~$\rho$ is the dual stepsize.
\begin{theorem}\label{thm:final}
Let Assumptions~\ref{as:f}-\ref{as:dualcomm} hold and fix~$\delta > 0$. 
Choose~$0 < \rho <  \frac{2 \delta}{\delta^2 + 2}$. Let~$T(t) := \min_c t_c$ be the minimum number of updates any one dual agent has performed by time~$t$ and let~$K(t)$ be the minimum value of~$\textnormal{ops}$ that was reached for 
any primal block used to compute any dual block from~$\mu(0)$ to~$\mu(t)$. 
Then for agents executing Algorithm~\ref{alg:2}, 
for all~$i$, all~$k$, and all~$t$, we have
\begin{equation}
\|x^i(k; t)\!-\!\hat{x}_{\delta}\|^2\! \leq\!q_p^{2 \textnormal{ops}(k,t)} 2 n D_x^2 
     + q_d^{T(t)} \frac{2 M^2} {\beta^2} \|\mu (0) - \mhatd \|^2  
     +  q_p^{2K(t)} C_1  + q_p^{K(t)} C_2  +  C_3,
\end{equation}
where~$C_1,$~$C_2$, and~$C_3$ are positive constants given by
\begin{equation}
    C_1 := \frac{2 n N_d M^4 D_{x}^2 (q_d - \rho^2)} {\beta^2 (1-q_d)}, \,\,\, C_2 := \frac{4 \rho^2 \sqrt{n} N_d M^4 D_{x}^2} {\beta^2 (1-q_d)}, \,\,\,
    C_3 := \frac{2 N_d  M^4 D_{x}^2 (q_d - \rho^2)} {\beta^2 (1-q_d)},
\end{equation}
and~$q_d := (1\!-\!\rho \delta)^2\!+\!2\rho^2 \in [0,1)$,~$q_p:=(1-\gamma \beta) \in [0,1)$,~$M := \max\limits_{x \in X} \| \nabla g(x) \|$,~$D_{x} := \max\limits_{x,y \in X} \| x-y \|$,~$n$ is the length of the primal variable~$x$, and~$N_d$ is the number of dual agents.
\end{theorem}}{\blue{We now present our main result on the convergence of~$x^i(k; t)$ to~$\xhat$. Recall that~$\delta$ is the dual regularization parameter,~$\gamma$ is the primal 
stepsize given in~\eqref{eq:gamma}, and~$\rho$ is the dual stepsize.
\begin{theorem}\label{thm:final}
Let Assumptions~\ref{as:f}-\ref{as:dualcomm} hold and fix~$\delta > 0$. 
Choose~$0 < \rho <  \frac{2 \delta}{\delta^2 + 2}$. Let~$T(t) := \min_c t_c$ be the minimum number of updates any one dual agent has performed by time~$t$ and let~$K(t)$ be the minimum value of~$\textnormal{ops}$ that was reached for 
any primal block used to compute any dual block from~$\mu(0)$ to~$\mu(t)$. 
Then for agents executing Algorithm~\ref{alg:2}, 
for all~$i$, all~$k$, and all~$t$, we have
\begin{multline}
\|x^i(k; t)\!-\!\hat{x}_{\delta}\|^2\! \leq\!q_p^{2 \textnormal{ops}(k,t)} 2 n D_x^2 
     + q_d^{T(t)} \frac{2 M^2} {\beta^2} \|\mu (0) - \mhatd \|^2  \\
     +  q_p^{2K(t)} C_1  + q_p^{K(t)} C_2  +  C_3,
\end{multline}
where~$C_1,$~$C_2$, and~$C_3$ are positive constants given by
\begin{align}
    C_1 &:= \frac{2 n N_d M^4 D_{x}^2 (q_d - \rho^2)} {\beta^2 (1-q_d)}, \,\,\, C_2 := \frac{4 \rho^2 \sqrt{n} N_d M^4 D_{x}^2} {\beta^2 (1-q_d)}, \\
    C_3 &:= \frac{2 N_d  M^4 D_{x}^2 (q_d - \rho^2)} {\beta^2 (1-q_d)},
\end{align}
and~$q_d := (1\!-\!\rho \delta)^2\!+\!2\rho^2 \in [0,1)$,~$q_p:=(1-\gamma \beta) \in [0,1)$,~$M := \max\limits_{x \in X} \| \nabla g(x) \|$,~$D_{x} := \max\limits_{x,y \in X} \| x-y \|$,~$n$ is the length of the primal variable~$x$, and~$N_d$ is the number of dual agents.
\end{theorem}}}
\blue{\emph{Proof:} See Appendix~\ref{app:overallconv}. \hfill $\blacksquare$}
\begin{remark}
\blue{The term~$C_3$ in Theorem~\ref{thm:final} is termed the ``asynchrony penalty''} because
it is an offset from reaching a solution, and it is not reduced by
changing the value of~$\textnormal{ops}$. 
It is due to asynchronously computing dual
variables and is absent when dual updates are centralized~\cite{koshal2011multiuser,hale17}. 
In Corollary~\ref{cor:asynchpen} below, we suggest methods to mitigate this penalty. 
\end{remark}

\blue{Note that the terms premultiplied by~$q_p$ are minimized when the exponent (determined by primal operations through the~$\textnormal{ops}$ term) is allowed to grow. In particular, if more primal operations occur before communications are sent to dual agents, the terms are reduced. Similarly, if dual updates occur frequently and asynchronously, the exponent will be reduced
and hence the terms containing~$q_p$ will become larger.}

\subsection{\blue{Reducing the Asynchrony Penalty}}
\blue{The asynchrony penalty may be addressed in a few ways: by completing a certain number of primal operations prior to sending updates to dual agents, by completing more dual updates, and finally by choosing~$\rho$ and~$\delta$. These are discussed explicitly in Corollary~\ref{cor:asynchpen} below.
\begin{corollary} \label{cor:asynchpen}
Let all conditions and definitions of Theorem~\ref{thm:final} hold. Let positive error bounds~$\epsilon_1$ and~$\epsilon_2$ be fixed. Then there exist values of~$T(t)$,~$K(t)$,
the dual stepsize~$\rho$, and the regularization parameter~$\delta$ such that
\begin{align}
   \|x^i(k; t)\!-\!\hat{x}_{\delta}\|^2  &\leq \epsilon_1 + \epsilon_2.
\end{align}
In particular, set~$K(t)\!\geq\!\frac{\log(\epsilon_1)\!-\!\log(4n D_x^2 \!+\!2C_1\!+\!2C_2)}{\log(q_p)}$,~$T(t) \geq \frac{\log (\epsilon_1 \beta^2) - \log(4 M^2 \|\mu (0) - \mhatd \|^2)}{\log (q_d)}$,~$\rho = \frac{\delta}{1+\delta^2}$, and~$\delta^2 \geq \frac{2 N_d  M^4 D_{x}^2}{\epsilon_2 \beta^2 (1-q_d) } - 1$.
\end{corollary}
\emph{Proof:} See Appendix~\ref{app:asynchpen}. \hfill $\blacksquare$}

\blue{The lower bound on~$\delta$ is dependent on the chosen~$\epsilon_2$. If~$\epsilon_2 \geq \frac{2 N_d  M^4 D_{x}^2} {\beta^2 (1-q_d)}$, then~$\delta^2$ may take any positive value. However, if~$\epsilon_2 \leq \frac{2 N_d  M^4 D_{x}^2} {\beta^2 (1-q_d)}$,~$\delta^2$ is lower bounded by a positive number. This illustrates the potential trade-off between the asynchrony penalty and regularization error by showing that
one cannot have arbitrarily small values of both. 
However, given a desired value of one, it is possible to compute the other
using Theorem~\ref{lem:regerror}, and,
if needed, make further adjustments to the numbers of updates and choices
of parameters to ensure satisfactory performance. }

\section{Numerical Example} \label{sec:numerical}
\begin{figure}[t!]
\centering
\includegraphics[width=5cm]{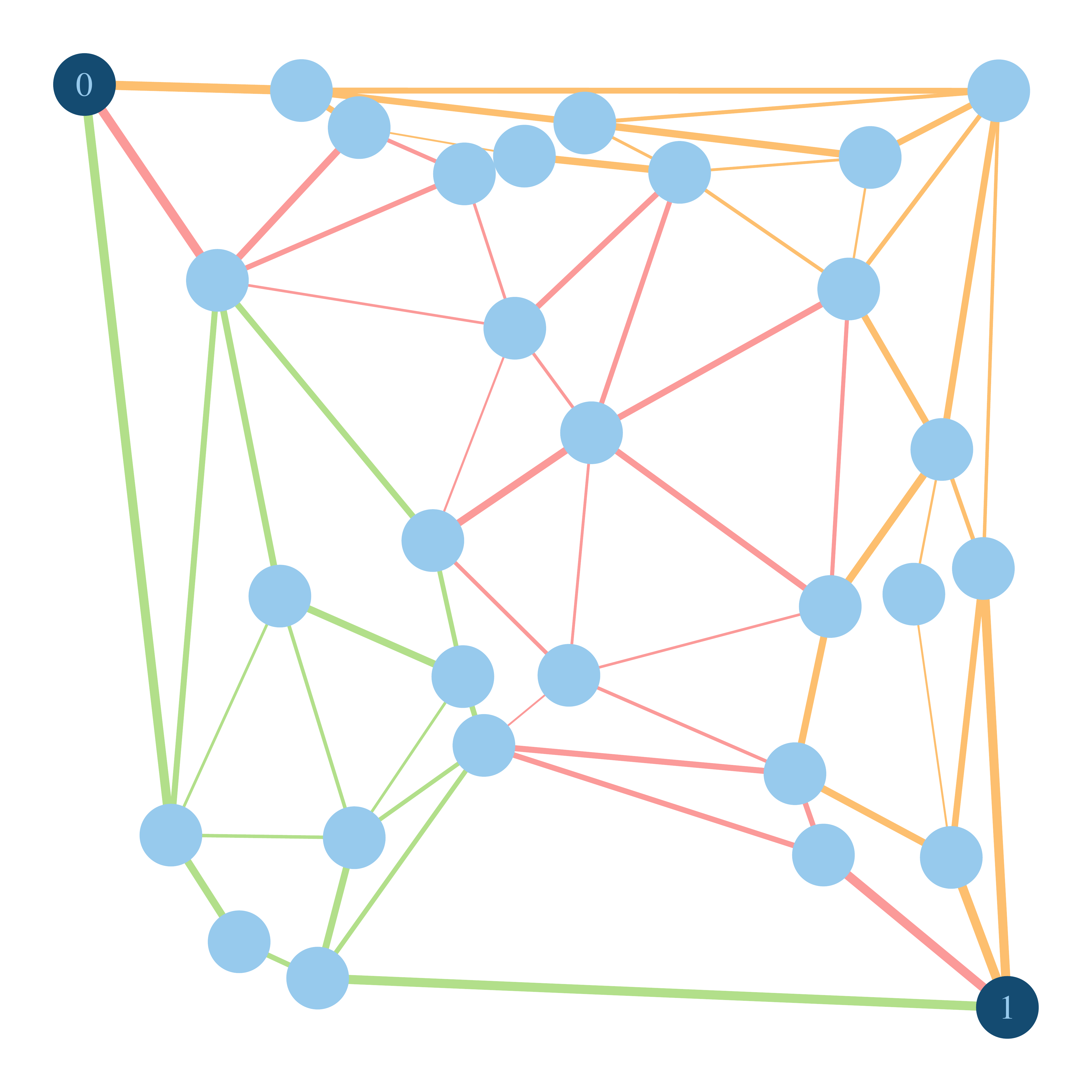}
\caption{Network graph where Node 0 is the source and Node 1 is the target and edge thicknesses correspond to their flow capacities. The paths may be divided into three groups (green, pink, and orange) such that there are no shared edges for paths in different groups.}
\label{fig:network}
\end{figure}

We consider a network flow problem where agents are attempting to route data over given paths from the source (node 0) to the target (node 1)\footnote{The network graph in Figure~\ref{fig:network} was generated using~\cite{graph-tool}. All other code for this section can be found at www.github.com/kathendrickson/dacoa.}. We consider an example whose network graph is given in Figure~\ref{fig:network}, where the edge widths represent the edge capacities. We consider a problem in which we must route~$15$ 
different paths from the source to the target along 
some combination of the 66 edges present. The primal variable is composed of the traffic assigned to each path (with~$n=15$) and the limits for traffic sent along each edge are the constraints (thus,~$m=66$). By construction, the paths and edges can be divided into three groups, with the edges in each group being used only by that group's paths. 
Each edge group (indicated by the green, pink, and orange colors in Figure~\ref{fig:network}) contains five paths, corresponding to five entries of the primal variable.  The objective function is
\begin{equation}
f(x) = - \blue{W} \sum_{i=1}^{n} \log (1+x_i),
\end{equation}
\blue{where~$W$ is a positive constant.}
Primal variables~$x_i$ are allowed to take any value between 0 and 10. 
The constraints are given by~$Ax \leq b$, where the edge capacities,~$b$, take random values between 5 and 40, with edges connected to the source and target having capacities of 50. The matrix~$A$ is given by
\begin{equation}
    A_{k,i} = \begin{cases}
    1 & \textnormal{if flow path }i \textnormal{ traverses edge }k \\
    0 & \textnormal{otherwise}
    \end{cases}.
\end{equation}
Thus the regularized Lagrangian is given by 
\begin{equation*}
L_{\delta}(x,\mu)
=  - \blue{W} \sum_{i=1}^{n} \log (1+x_i) + \mu^T\Big(Ax-b\Big) - \frac{\delta}{2}\|\mu\|^2, 
\end{equation*}
where the Hessian matrix~$H$ is~$\beta$-diagonally dominant \blue{and~$\beta = \frac{W}{11^2}$.}
We choose algorithm parameters~$\gamma = 0.01$,~$\delta = 0.1$ and~$\rho = \blue{ \frac{\delta}{\delta^2 + 1} \approx .099}$. Communications between agents occur with a random probability called the ``communication rate,'' 
which we vary across simulation runs below. \blue{We collect each primal agent's block into the combined primal variable~$x = (x^{1T}_{[1]}, \dots, x^{N_pT}_{[N_p]})^T$ to measure convergence.}

We use this simulation example to explore the benefit of using non-scalar blocks over our previous work with scalar blocks~\cite{hendrickson20}. Additionally, we examine the effect that the magnitude of diagonal dominance has on convergence. We are also able to vary the communication rate and study its effect on convergence. Finally, this example demonstrates the effectiveness of our algorithm with large-scale problems and the ease with which it is scaled up and distributed.

We begin by comparing scalar blocks to non-scalar blocks \blue{where primal agents have a 50\% chance of computing an update at every time~$k$}. When using scalar blocks, we assign one primal agent to each flow path and one dual agent to each edge constraint. Thus, we have 15 primal agents and 66 dual agents in the scalar block case. For dividing among non-scalar blocks,  we assign one primal agent to compute all five flow paths in each network group (indicated by the different colors in Figure~\ref{fig:network}). We then assign one dual agent to each group to handle all of the edge constraints for that group. Thus, we have 3 primal agents and 3 dual agents in the non-scalar block case.  For a communication rate of 0.75 (primal agents have a 75\% chance of communicating the latest update to another agent at each \blue{time step}), non-scalar blocks provide an advantage when considering the  number of \blue{time steps} needed to converge, shown in Figure~\ref{fig:blocks}. \blue{In both cases, the algorithm converges to~$\tilde{x}$ such that~$\|\tilde{x} - \hat{x}\| \leq 0.38$, where~$\hat{x}$ is the unregularized solution. This result is representative of other simulations by the authors in which the asynchrony penalty is small in practice and the bound provided in Theorem~\ref{thm:final} is loose.}

Next we \blue{use the non-scalar blocks with primal agents performing updates at every time~$k$ to isolate the effects of diagonal dominance and communication rate (communications do not necessarily happen at every~$k$). We first} vary~$\beta$ over\blue{~$\beta \in \lbrace 0.10, 0.25, 0.75\rbrace$.} \blue{To measure convergence, we take the~$2$-norm of the difference between~$x$ at consecutive time steps. Figure~\ref{fig:beta} plots the \blue{time step}~$k$ versus this distance for the varying values of~$\beta$ and a communication rate of~$0.75$}. As predicted by \blue{Theorem~\ref{thm:final}}, a larger~$\beta$ correlates with faster convergence in general. 
%\red{However, as~$\beta$ grows so do oscillations. Depending on the context, the number of iterations must be weighed with the effect of oscillations during implementation. For purely analytical applications, oscillations have little impact but reducing the number of iterations is highly desirable. However, in physical applications such as robotics, oscillations may create erratic behavior that is not worth reduced iterations.} 

\ifbool{Report}{\begin{figure}[t!]
\centering
\includegraphics[width=8.8cm]{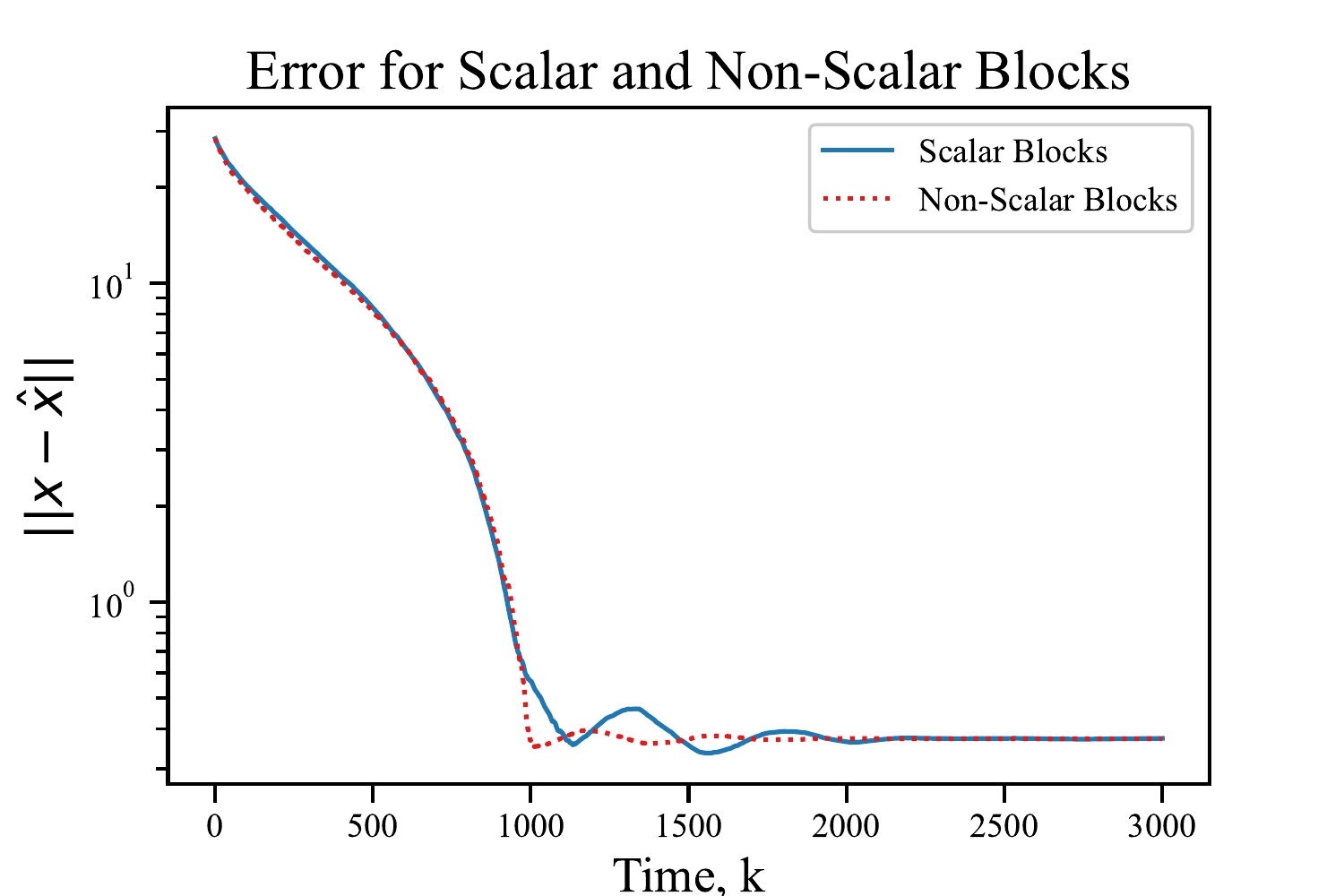}
\caption{Convergence for scalar and non-scalar blocks. Non-scalar blocks provide a significant advantage over scalar blocks when considering the number of \blue{time steps} needed to reach a solution, as indicated by the red dotted line versus the solid blue line.}
\label{fig:blocks}
\end{figure}
\begin{figure}[t!]
\centering
\includegraphics[width=8.8cm]{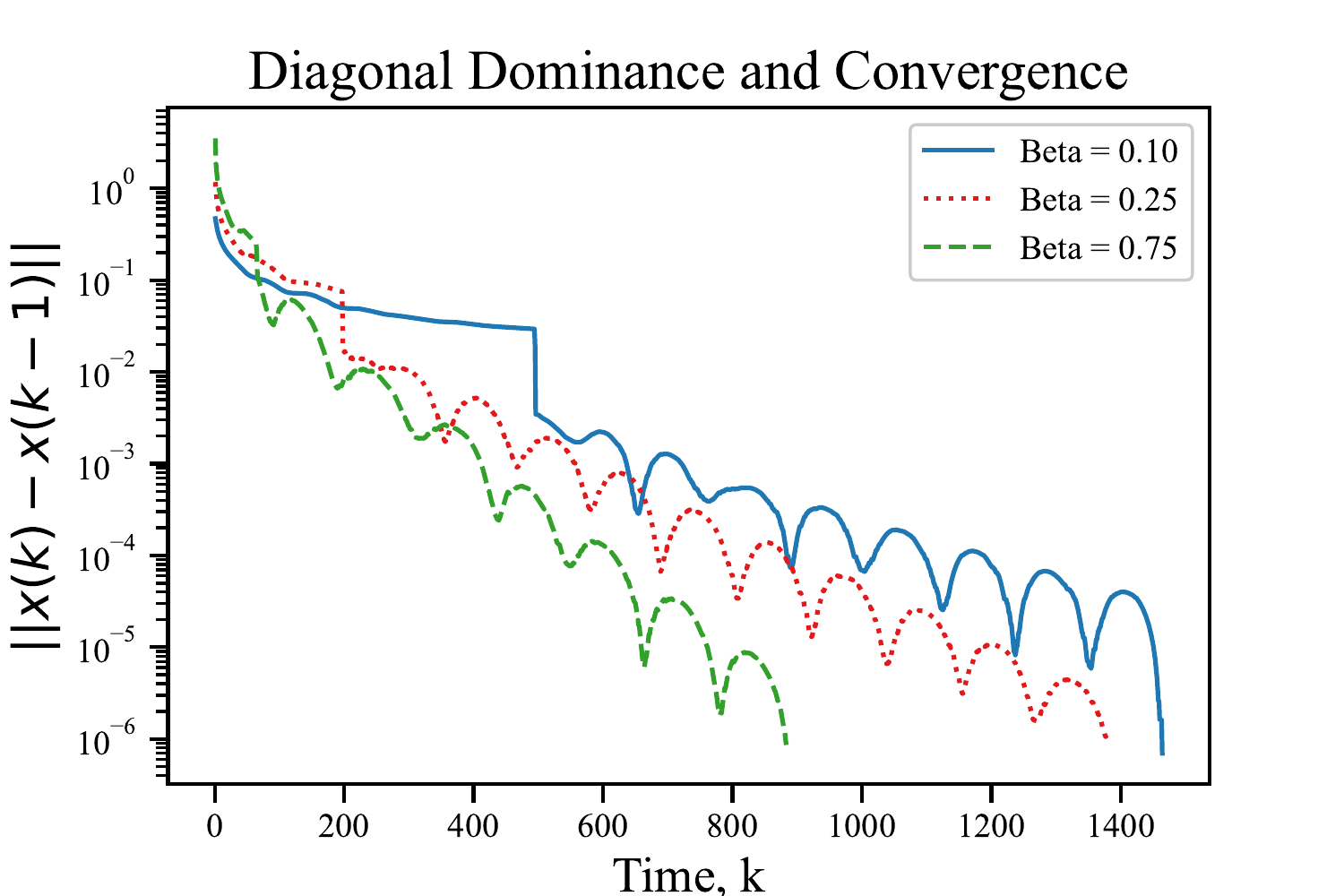}
\caption{Effect of diagonal dominance on convergence. Here, we see that larger values of~$\beta$ lead to faster convergence.}
\label{fig:beta}
\end{figure}
\begin{figure}[t!]
\centering
\includegraphics[width=8.8cm]{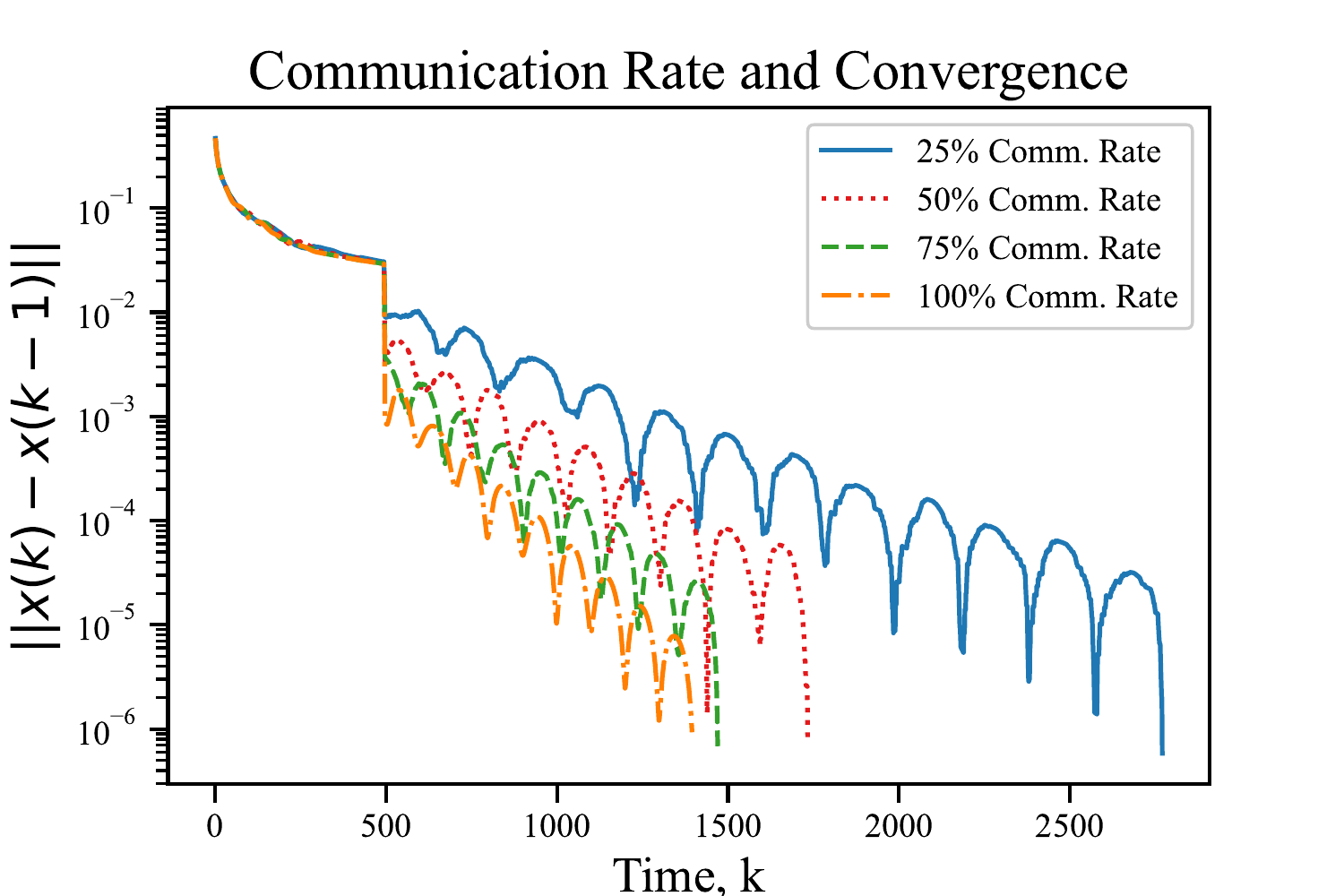}
\caption{Effect of communication rate on convergence. Less frequent communication leads to slower convergence.}
\label{fig:comm}
\end{figure}}{\begin{figure}[t!]
\centering
\includegraphics[width=8.8cm]{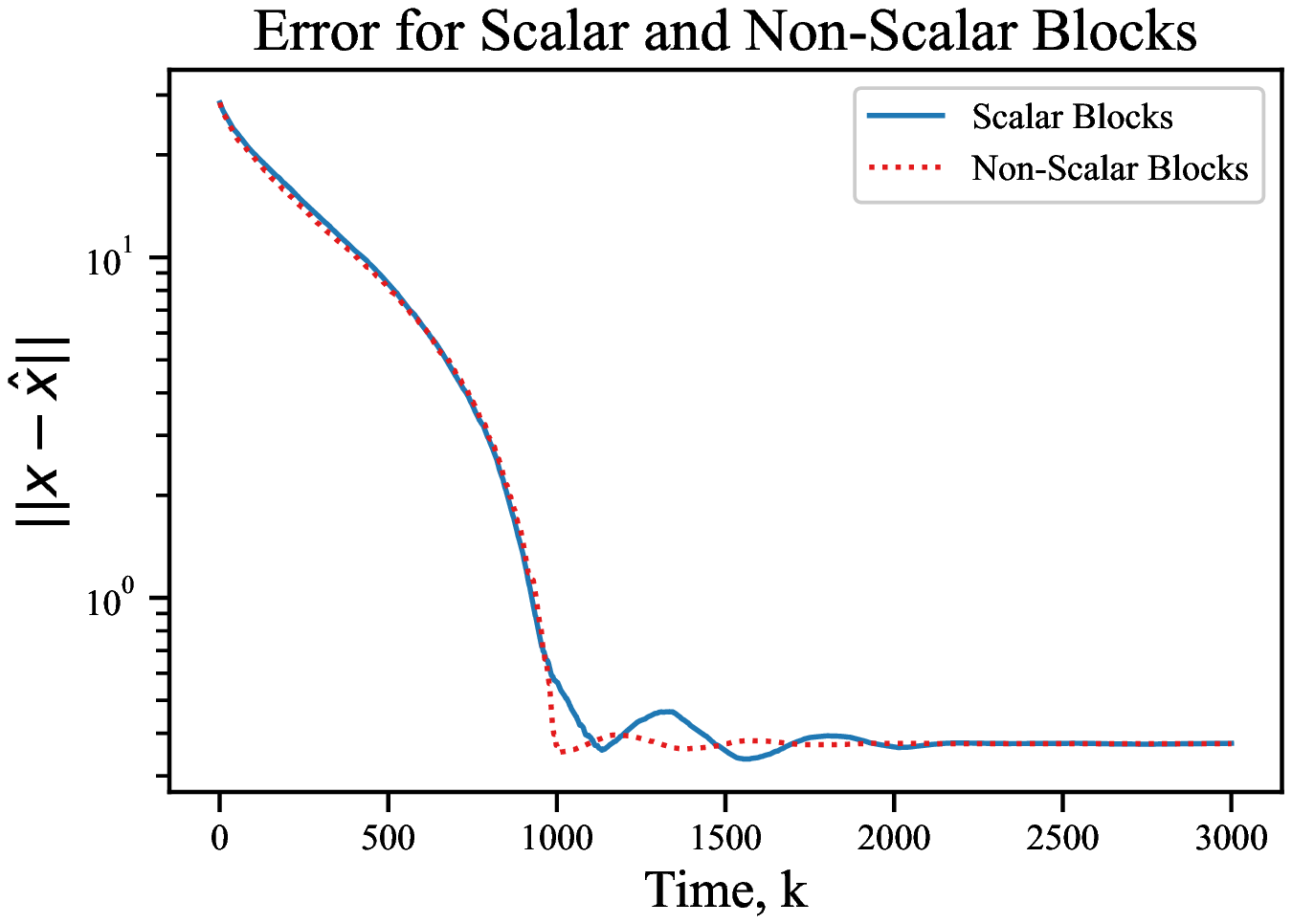}
\caption{Convergence for scalar and non-scalar blocks. Non-scalar blocks provide a significant advantage over scalar blocks when considering the number of \blue{time steps} needed to reach a solution, as indicated by the red dotted line versus the solid blue line.}
\label{fig:blocks}
\end{figure}
\begin{figure}[t!]
\centering
\includegraphics[width=8.8cm]{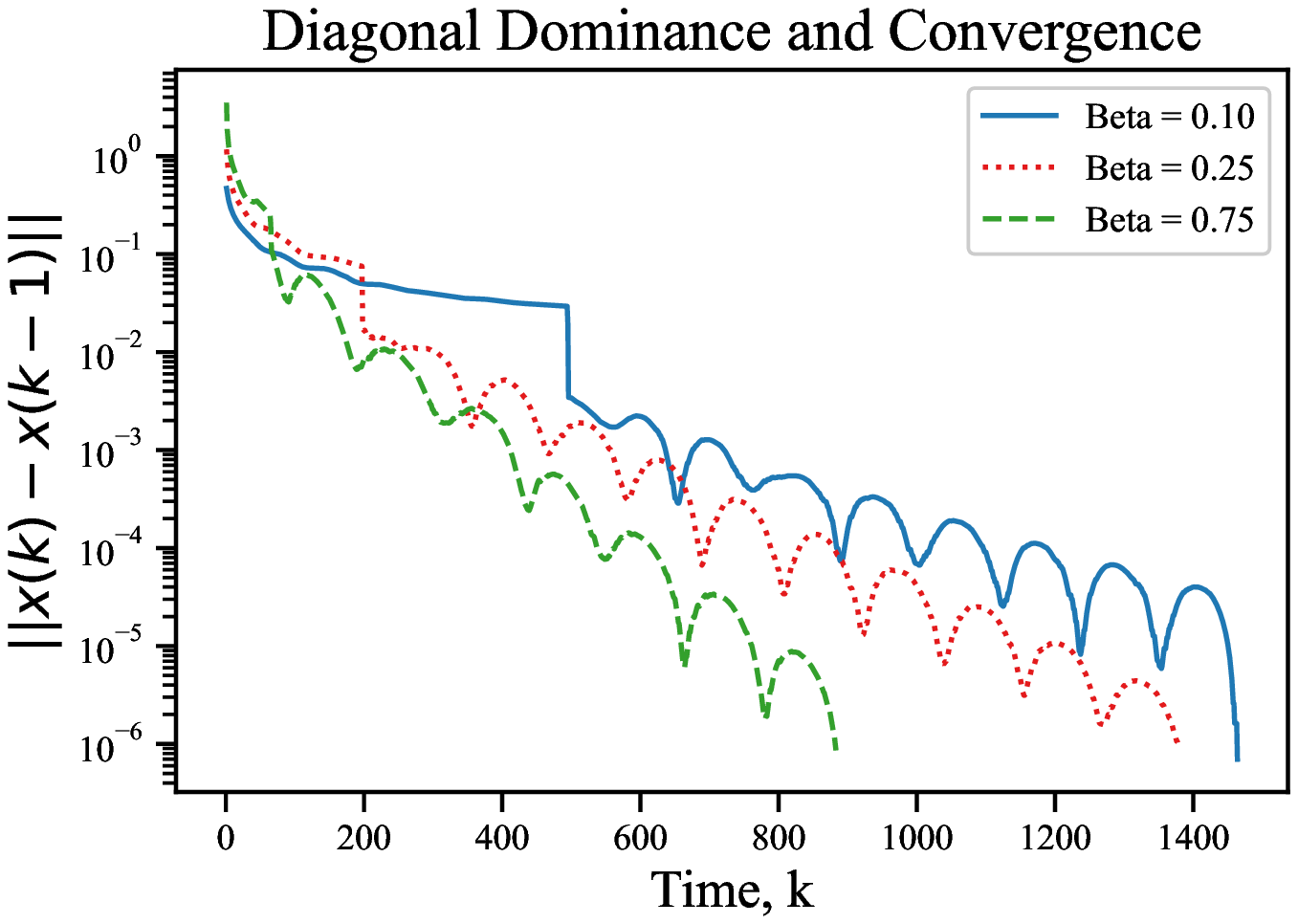}
\caption{Effect of diagonal dominance on convergence. Here, we see that larger values of~$\beta$ lead to faster convergence.}
\label{fig:beta}
\end{figure}
\begin{figure}[t!]
\centering
\includegraphics[width=8.8cm]{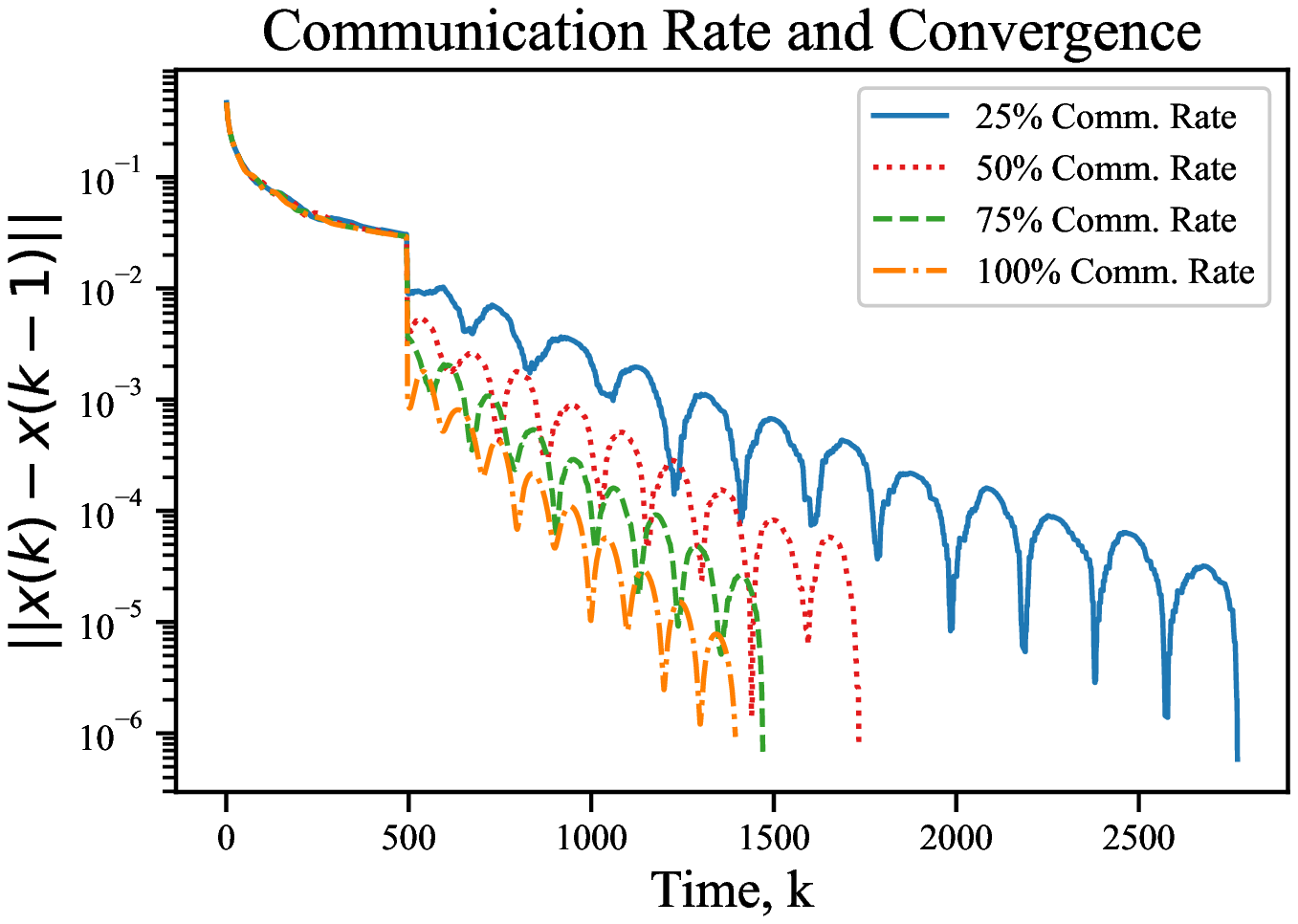}
\caption{Effect of communication rate on convergence. Less frequent communication leads to slower convergence.}
\label{fig:comm}
\end{figure}}

Varying the communication rate has a significant impact on the number of \blue{time steps required to converge as shown in Figure~\ref{fig:comm}. However, a solution is still eventually reached.} This reveals that faster convergence can be achieved both by increasing communication rates and increasing the diagonal dominance of the problem. 
\blue{
In these two plots, the abrupt decreases in distances between successive iterates are due to primal agents' computations reaching the boundary
of the feasible region defined by~$X$ and~$g$, which causes the iterates to make only small progress
afterwards (as the dual variables continue to slowly change). 
}

\section{Conclusion} \label{sec:concl}
\blue{Algorithm~\ref{alg:2} presents a primal-dual approach that is asynchronous in primal updates and communications and asynchronous in distributed dual updates. The error due to regularization was bounded and convergence rates were established. Methods for mitigating the resulting asynchrony penalty were presented. A numerical example illustrates the benefit of non-scalar blocks and the effect diagonal dominance has with other parameters upon convergence. Future work will examine additional applications for the algorithm and implementation techniques to reduce error and improve convergence.}

\appendix
\subsection{Proof of Theorem~\ref{lem:regerror}} \label{app:regerror}
\ifbool{Report}{The following proof generally follows that of Proposition 3.1 in~\cite{koshal2011multiuser}, with differences resulting from this work
only 
regularizing with respect to the dual variable rather than regularizing with respect to both the primal and dual variables.
\blue{Let~$(\hat{x}, \hat{\mu})$ denote a saddle point of the unregularized Lagrangian~$L$. Let~$(\hat{x}_\delta, \hat{\mu}_\delta)$ denote a saddle point of the dual-regularized Lagrangian~$L_\delta$. Then because~$(\hat{x}_\delta,\hat{\mu}_\delta)$ is a saddle point, for all~$x\in X, \mu \in \R^m_+$ we have the two inequalities
%\begin{align}
    $L_\delta (\hat{x}_\delta, \mu) \leq L_\delta (\hat{x}_\delta, \hat{\mu}_\delta) \leq L_\delta (x, \hat{\mu}_\delta)$. % \label{eq:saddle}
%\end{align}
Using~$\hat{\mu} \in \R^m_+$ and we can write
\begin{align}
    0 &\leq L_\delta (\hat{x}_\delta, \hat{\mu}_\delta) - L_\delta (\hat{x}_\delta, \hat{\mu}) =\sum_j(\hat{\mu}_{\delta,j} - \hat{\mu}_j) g_j(\hat{x}_\delta) - \frac{\delta}{2}\|\hat{\mu}_\delta\|^2 + \frac{\delta}{2}\|\hat{\mu}\|^2. \label{eq:saddleexp}
\end{align}
Because each~$g_j$ is convex, we have
\begin{equation}
    g_j(\hat{x}_\delta) \!\leq g_j(\hat{x}) \!+\! \nabla g_j(\hat{x}_\delta)^T \!(\hat{x}_\delta - \hat{x}) \!\leq \!\nabla g_j(\hat{x}_\delta)^T (\hat{x}_\delta - \hat{x}), \label{eq:gconvex}
\end{equation}
where the last inequality follows from~$g_j(\hat{x}) \leq 0$ (which holds
since~$\hat{x}$ solves Problem~\ref{prob:first}). 
%as~$\hat{x}$ is the solution to the original Lagrangian (and thus, the constraints are satisfied). 
Additionally, because all dual variables are non-negative, 
we can multiply by~$\hat{\mu}_{\delta,j}$ to get
\begin{equation}
    \sum_j \hat{\mu}_{\delta,j}  g_j(\hat{x}_\delta) \leq \sum_j \hat{\mu}_{\delta,j} \nabla g_j(\hat{x}_\delta)^T (\hat{x}_\delta - \hat{x}). \label{eq:appAdualineqbound}
\end{equation}
By definition of~$L_\delta$, the right-hand side can be expanded as
\begin{align}
    \sum_j \hat{\mu}_{\delta,j} \nabla g_j(\hat{x}_\delta)^T (\hat{x}_\delta - \hat{x}) &= \nabla_x L_\delta(\hat{x}_\delta, \hat{\mu}_{\delta})^T (\hat{x}_\delta - \hat{x}) - \nabla f(\hat{x}_\delta)^T (\hat{x}_\delta - \hat{x}). \label{eq:appAgrad}
\end{align}
Because~$\xhat$ minimizes~$L_{\delta}(\cdot, \mhatd)$, for all~$x \in X$
we have~$\nabla_x L_\delta(\hat{x}_\delta, \hat{\mu}_{\delta})^T (\hat{x}_\delta - x) \leq 0$.
In particular, we can set~$x = \hat{x} \in X$ to find~$\nabla_x L_\delta(\hat{x}_\delta, \hat{\mu}_{\delta})^T (\hat{x}_\delta - \hat{x}) \leq 0$. 
Combining this with~\eqref{eq:appAgrad} and~\eqref{eq:appAdualineqbound} gives 
\begin{align} \label{eq:sum1} 
    \sum_j \hat{\mu}_{\delta,j}  g_j(\hat{x}_\delta) \leq - \nabla f(\hat{x}_\delta)^T (\hat{x}_\delta - \hat{x}). 
\end{align}
By the convexity of each~$g_j$, we have~$g_j(\hat{x}_\delta) \geq g_j(\hat{x}) + \nabla g_j(\hat{x})^T(\hat{x}_\delta - \hat{x})$. By multiplying this inequality with the non-positive~$-\hat{\mu}_j$ and summing over~$j$, we obtain
\begin{equation} \label{eq:forL1} 
    -\!\sum_j \hat{\mu}_j g_j(\hat{x}_\delta) \leq\!-\!\sum_j \hat{\mu}_j g_j(\hat{x})\!-\!\sum_j \hat{\mu}_j \nabla g_j(\hat{x})^T(\hat{x}_\delta\!-\!\hat{x}).
\end{equation}
By complementary slackness, we have~$\hat{\mu}^T g(\hat{x}) = 0$ and thus 
\begin{equation} \label{eq:Lsum2}
    -\sum_j \hat{\mu}_j g_j(\hat{x}_\delta) \leq - \sum_j \hat{\mu}_j \nabla g_j(\hat{x})^T(\hat{x}_\delta - \hat{x}).
\end{equation}
Expanding~$\nabla_x L(\hat{x}, \hat{\mu}) = \nabla f(\hat{x}) + \sum_j \hat{\mu}_j \nabla g_j(\hat{x})$, we see that
\begin{equation}
    \sum_j \hat{\mu}_j \nabla g_j(\hat{x})^T  (\hat{x} - \hat{x}_\delta) =   \nabla_x L(\hat{x}, \hat{\mu})^T (\hat{x} - \hat{x}_\delta) - \nabla  f(\hat{x})^T (\hat{x} - \hat{x}_\delta) 
    \leq - \nabla f(\hat{x})^T(\hat{x} - \hat{x}_\delta).
\end{equation}
This follows from the fact that~$\hat{x}$ minimizes~$L(\cdot, \hat{\mu})$ over all~$x \in X$
and thus~$\nabla_x L(\hat{x}, \hat{\mu})^T(\hat{x} - x) \leq 0$ for all~$x \in X$.
Then setting~$x = \xhat$ gives the above bound. 
Combining this with~\eqref{eq:Lsum2} gives 
\begin{equation}
    -\sum_j \hat{\mu}_j g_j(\hat{x}_\delta) %&\leq - \nabla f(\hat{x})^T(\hat{x} - \hat{x}_\delta) 
    \leq \nabla f(\hat{x})^T(\hat{x}_\delta - \hat{x}). \label{eq:sum2}
\end{equation}
Adding~\eqref{eq:sum1} and~\eqref{eq:sum2} gives
\begin{equation}
    %&\sum_j \hat{\mu}_{\delta,j}  g_j(\hat{x}_\delta) -\sum_j \hat{\mu}_j g_j(\hat{x}_\delta) \leq  - \nabla_x f(\hat{x}_\delta)^T (\hat{x}_\delta - \hat{x}) \\
    %&\qquad \qquad \qquad + \nabla_x f(\hat{x})^T(\hat{x}_\delta - \hat{x}) \\
    (\hat{\mu}_{\delta} - \hat{\mu})^T g(\hat{x}_\delta) \leq  \left(\nabla f(\hat{x}) - \nabla f(\hat{x}_\delta) \right)^T (\hat{x}_\delta - \hat{x}) 
     \leq -\frac{\beta}{2}\|\hat{x}_\delta - \hat{x}\|^2,
\end{equation}
where the last inequality is from the~$\beta$-strong convexity of~$f$ (follows from Assumption~\ref{as:diagonal} which holds for~$\mu=0$). 
Applying this to~\eqref{eq:saddleexp},
\begin{align}
    0 &\leq -\frac{\beta}{2}\|\hat{x}_\delta - \hat{x}\|^2 - \frac{\delta}{2}\|\hat{\mu}_\delta\|^2 + \frac{\delta}{2}\|\hat{\mu}\|^2.
\end{align}
This implies the final result~$\|\hat{x}_\delta - \hat{x}\|^2 \leq \frac{\delta}{\beta}\left(\|\hat{\mu}\|^2 - \|\hat{\mu}_\delta\|^2\right).$
% \begin{align}
%     \|\hat{x}_\delta - \hat{x}\|^2 \leq \frac{\delta}{\beta}\left(\|\hat{\mu}\|^2 - \|\hat{\mu}_\delta\|^2\right).
% \end{align}
Furthermore, we can use~\eqref{eq:gconvex} to bound possible constraint violations, where
\begin{align}
    g_j(\hat{x}_\delta) &\leq \nabla g_j(\hat{x}_\delta)^T (\hat{x}_\delta - \hat{x}) \leq \|\nabla g_j(\hat{x}_\delta)\| \|\hat{x}_\delta - \hat{x}\| \leq \max_{x \in X} \|\nabla g_j(x)\| \sqrt{\frac{\delta  }{\beta}\left(\|\hat{\mu}\|^2 - \|\hat{\mu}_\delta\|^2\right)} \leq M_j \munormbound\sqrt{\frac{\delta  }{\beta}}. \tag*{$\blacksquare$} 
\end{align}}}{\blue{The following proof generally follows that of Proposition 3.1 in~\cite{koshal2011multiuser}, with differences resulting from this work
only 
regularizing with respect to the dual variable rather than regularizing with respect to both the primal and dual variables.}
\blue{Let~$(\hat{x}, \hat{\mu})$ denote a saddle point of the unregularized Lagrangian~$L$. Let~$(\hat{x}_\delta, \hat{\mu}_\delta)$ denote a saddle point of the dual-regularized Lagrangian~$L_\delta$. Then because~$(\hat{x}_\delta,\hat{\mu}_\delta)$ is a saddle point, for all~$x\in X, \mu \in \R^m_+$ we have the two inequalities
%\begin{align}
    $L_\delta (\hat{x}_\delta, \mu) \leq L_\delta (\hat{x}_\delta, \hat{\mu}_\delta) \leq L_\delta (x, \hat{\mu}_\delta)$. % \label{eq:saddle}
%\end{align}
Using~$\hat{\mu} \in \R^m_+$ and we can write
\begin{align}
    0 &\leq L_\delta (\hat{x}_\delta, \hat{\mu}_\delta) - L_\delta (\hat{x}_\delta, \hat{\mu}) \\
    %&= (\hat{\mu}_\delta - \hat{\mu})^T g(\hat{x}_\delta) - \frac{\delta}{2}\|\hat{\mu}_\delta\|^2 + \frac{\delta}{2}\|\hat{\mu}\|^2 \\
    &=\sum_j(\hat{\mu}_{\delta,j} - \hat{\mu}_j) g_j(\hat{x}_\delta) - \frac{\delta}{2}\|\hat{\mu}_\delta\|^2 + \frac{\delta}{2}\|\hat{\mu}\|^2. \label{eq:saddleexp}
\end{align}
Because each~$g_j$ is convex, we have
\begin{equation}
    g_j(\hat{x}_\delta) \!\leq g_j(\hat{x}) \!+\! \nabla g_j(\hat{x}_\delta)^T \!(\hat{x}_\delta - \hat{x}) \!\leq \!\nabla g_j(\hat{x}_\delta)^T (\hat{x}_\delta - \hat{x}), \label{eq:gconvex}
\end{equation}
where the last inequality follows from~$g_j(\hat{x}) \leq 0$ (which holds
since~$\hat{x}$ solves Problem~\ref{prob:first}). 
%as~$\hat{x}$ is the solution to the original Lagrangian (and thus, the constraints are satisfied). 
Additionally, because all dual variables are non-negative, 
we can multiply by~$\hat{\mu}_{\delta,j}$ to get
\begin{equation}
    \sum_j \hat{\mu}_{\delta,j}  g_j(\hat{x}_\delta) \leq \sum_j \hat{\mu}_{\delta,j} \nabla g_j(\hat{x}_\delta)^T (\hat{x}_\delta - \hat{x}). \label{eq:appAdualineqbound}
\end{equation}
By definition of~$L_\delta$, the right-hand side  can be expanded as
\begin{align}
    \sum_j \hat{\mu}_{\delta,j} \nabla g_j(\hat{x}_\delta)^T (\hat{x}_\delta - \hat{x}) &= \nabla_x L_\delta(\hat{x}_\delta, \hat{\mu}_{\delta})^T (\hat{x}_\delta - \hat{x}) \\
    &\qquad - \nabla f(\hat{x}_\delta)^T (\hat{x}_\delta - \hat{x}). \label{eq:appAgrad}
\end{align}
Because~$\xhat$ minimizes~$L_{\delta}(\cdot, \mhatd)$, for all~$x \in X$
we have~$\nabla_x L_\delta(\hat{x}_\delta, \hat{\mu}_{\delta})^T (\hat{x}_\delta - x) \leq 0$.
In particular, we can set~$x = \hat{x} \in X$ to find~$\nabla_x L_\delta(\hat{x}_\delta, \hat{\mu}_{\delta})^T (\hat{x}_\delta - \hat{x}) \leq 0$. 
Combining this with~\eqref{eq:appAgrad} and~\eqref{eq:appAdualineqbound} gives 
\begin{align} \label{eq:sum1} 
    \sum_j \hat{\mu}_{\delta,j}  g_j(\hat{x}_\delta) \leq - \nabla f(\hat{x}_\delta)^T (\hat{x}_\delta - \hat{x}). 
\end{align}
By the convexity of each~$g_j$, we have~$g_j(\hat{x}_\delta) \geq g_j(\hat{x}) + \nabla g_j(\hat{x})^T(\hat{x}_\delta - \hat{x})$. By multiplying this inequality with the non-positive~$-\hat{\mu}_j$ and summing over~$j$, we obtain
\begin{equation} \label{eq:forL1} 
    -\!\sum_j \hat{\mu}_j g_j(\hat{x}_\delta) \leq\!-\!\sum_j \hat{\mu}_j g_j(\hat{x})\!-\!\sum_j \hat{\mu}_j \nabla g_j(\hat{x})^T(\hat{x}_\delta\!-\!\hat{x}).
\end{equation}
By complementary slackness, we have~$\hat{\mu}^T g(\hat{x}) = 0$ and thus 
\begin{equation} \label{eq:Lsum2}
    -\sum_j \hat{\mu}_j g_j(\hat{x}_\delta) \leq - \sum_j \hat{\mu}_j \nabla g_j(\hat{x})^T(\hat{x}_\delta - \hat{x}).
\end{equation}
Expanding~$\nabla_x L(\hat{x}, \hat{\mu}) = \nabla f(\hat{x}) + \sum_j \hat{\mu}_j \nabla g_j(\hat{x})$, we see that
\begin{multline}
    \!\!\!\!\sum_j\!\hat{\mu}_j\!\nabla g_j(\hat{x})^T\! (\hat{x}\!-\!\hat{x}_\delta)\!=\!  \nabla_x L(\hat{x}, \hat{\mu})^T\!(\hat{x}\!-\!\hat{x}_\delta)\!-\!\nabla\! f(\hat{x})^T\!(\hat{x}\!-\!\hat{x}_\delta) \\
    \leq - \nabla f(\hat{x})^T(\hat{x} - \hat{x}_\delta).
\end{multline}
This follows from the fact that~$\hat{x}$ minimizes~$L(\cdot, \hat{\mu})$ over all~$x \in X$
and thus~$\nabla_x L(\hat{x}, \hat{\mu})^T(\hat{x} - x) \leq 0$ for all~$x \in X$.
Then setting~$x = \xhat$ gives the above bound. 
Combining this with~\eqref{eq:Lsum2} gives 
\begin{equation}
    -\sum_j \hat{\mu}_j g_j(\hat{x}_\delta) %&\leq - \nabla f(\hat{x})^T(\hat{x} - \hat{x}_\delta) 
    \leq \nabla f(\hat{x})^T(\hat{x}_\delta - \hat{x}). \label{eq:sum2}
\end{equation}
Adding~\eqref{eq:sum1} and~\eqref{eq:sum2} gives
\begin{multline}
    %&\sum_j \hat{\mu}_{\delta,j}  g_j(\hat{x}_\delta) -\sum_j \hat{\mu}_j g_j(\hat{x}_\delta) \leq  - \nabla_x f(\hat{x}_\delta)^T (\hat{x}_\delta - \hat{x}) \\
    %&\qquad \qquad \qquad + \nabla_x f(\hat{x})^T(\hat{x}_\delta - \hat{x}) \\
    \!\!\!\!\!(\hat{\mu}_{\delta} - \hat{\mu})^T g(\hat{x}_\delta)\!\leq\! \left(\nabla f(\hat{x})\!-\!\nabla f(\hat{x}_\delta) \right)^T (\hat{x}_\delta - \hat{x}) 
    \!\leq\!-\frac{\beta}{2}\|\hat{x}_\delta - \hat{x}\|^2,
\end{multline}
where the last inequality is from the~$\beta$-strong convexity of~$f$ (follows from Assumption~\ref{as:diagonal} which holds for~$\mu=0$). 
Applying this to~\eqref{eq:saddleexp},
\begin{align}
    0 &\leq -\frac{\beta}{2}\|\hat{x}_\delta - \hat{x}\|^2 - \frac{\delta}{2}\|\hat{\mu}_\delta\|^2 + \frac{\delta}{2}\|\hat{\mu}\|^2.
\end{align}
This implies the final result~$\|\hat{x}_\delta - \hat{x}\|^2 \leq \frac{\delta}{\beta}\left(\|\hat{\mu}\|^2 - \|\hat{\mu}_\delta\|^2\right).$
% \begin{align}
%     \|\hat{x}_\delta - \hat{x}\|^2 \leq \frac{\delta}{\beta}\left(\|\hat{\mu}\|^2 - \|\hat{\mu}_\delta\|^2\right).
% \end{align}
Furthermore, we can use~\eqref{eq:gconvex} to bound possible constraint violations, where
\ifbool{Report}{\begin{align}
    &g_j(\hat{x}_\delta) \leq \nabla g_j(\hat{x}_\delta)^T (\hat{x}_\delta - \hat{x}) \leq \|\nabla g_j(\hat{x}_\delta)\| \|\hat{x}_\delta - \hat{x}\| \\
    &\leq \max_{x \in X} \|\nabla g_j(x)\| \sqrt{\frac{\delta  }{\beta}\left(\|\hat{\mu}\|^2 - \|\hat{\mu}_\delta\|^2\right)} 
    \leq M_j \munormbound\sqrt{\frac{\delta  }{\beta}}. \tag*{$\blacksquare$} 
\end{align}}{\begin{align}
    &g_j(\hat{x}_\delta) \leq \nabla g_j(\hat{x}_\delta)^T (\hat{x}_\delta - \hat{x}) \leq \|\nabla g_j(\hat{x}_\delta)\| \|\hat{x}_\delta - \hat{x}\| \\
    &\leq \max_{x \in X} \|\nabla g_j(x)\| \sqrt{\frac{\delta  }{\beta}\left(\|\hat{\mu}\|^2 - \|\hat{\mu}_\delta\|^2\right)} 
    \leq M_j \munormbound\sqrt{\frac{\delta  }{\beta}}. \tag*{$\blacksquare$} 
\end{align}} 
}}

\subsection{Main Result Proofs} \label{app:overallconv}
\subsubsection{\blue{Primal Convergence}} \label{app:primal}
Towards defining an overall convergence rate to the regularized optimal solution~$(\xhat , \mhatd)$, we first find the primal convergence rate for a fixed dual variable. 
Given a fixed~$\mu(t)$, centralized projected gradient descent for minimizing~$L_{\delta}(\cdot,\mu(t))$ can be written as
\begin{equation} \label{eq:grad}
h(x) = \Pi_{X}\left[x - \gamma\nabla_{x}L_{\delta}(x,\mu(t))\right], 
\end{equation}
where~$\gamma > 0$. The fixed point of~$h$ is the minimizer of~$L_{\delta}(\cdot,\mu(t))$ and is denoted by~$\xhatt = \argmin_{x \in X} L_{\delta}(x, \mu(t))$.

Leveraging some existing theoretical tools in the study of
optimization algorithms~\cite{Bertsekas1983, Bertsekas1991}, we can
study~$h$ in a way that elucidates its behavior under asynchrony in a distributed implementation. 
%We do so using the max norm, defined for~$v \in \R^n$ as~$\|v\|_\infty = \max_i |v_i|$. 
According to \cite{Bertsekas1991}, the assumption of diagonal dominance guarantees that~$h$ has the contraction property \blue{
\begin{equation}
{\Vert h(x)-\xhatt \Vert_\infty \leq \alpha \Vert x - \xhatt \Vert}_\infty  
\end{equation}
for all~$x \in X$, 
where~$\Vert v\Vert _\infty = \max\limits_i |v_i|$ for~$v \in \R^n$,}~$\alpha \in [0,1)$ and~$\xhatt$ is a fixed point
of~$h$, which depends on the choice of fixed~$\mu(t)$.  
However, the value of~$\alpha$ is not specified in~\cite{Bertsekas1991}, 
and it is precisely that value that governs the rate of convergence
to a solution. We therefore compute~$\alpha$ explicitly. 

Following the method in \cite{Bertsekas1983}, two~$n \times n$ matrices~$G$ and~$F$ must also be defined.
\begin{definition} \label{def:GHmatrices} Define the~$n \times n$ matrices~$G$ and~$F$ as
\begin{equation}
   G=
    \begin{bmatrix}
    |H_{11}| & -|H_{12}| & \ldots & -|H_{1n}| \\
    \vdots & \vdots & \ddots & \vdots  \\
    -|H_{n1}| & -|H_{n2}| & \ldots & |H_{nn}| \\
    \end{bmatrix} \,\, \textnormal{and} \,\, F = I - \gamma G,
    %F = I\!-\!\gamma G &= \begin{bmatrix}
%1\!-\!\gamma |H_{11}| & \gamma|H_{12}| & \ldots & \gamma|H_{1n}| \\
%\vdots &\vdots &\ddots & \vdots \\
%\gamma|H_{n1}| & \gamma|H_{n2}| & \ldots & 1\!-\!\gamma |H_{nn}|,
%\end{bmatrix}
\end{equation}
where~$I$ is the~$n \times n$ identity matrix.
\end{definition}

We now have the following. 

\begin{lemma} \label{83lemma}
Let~$h$,~$G$, and~$F$ be as above and let \blue{Assumptions~\ref{as:f}-\ref{as:dualcomm}} hold. 
Then
%\begin{equation}
$|h(x) - h(y)| \leq F|x-y|$
%\end{equation}
 for all~$x,y \in \R^n$, 
where~$|v|$ denotes the element-wise absolute value
of the vector~$v \in \mathbb{R}^n$ and the inequality holds component-wise. 
\end{lemma}
\emph{Proof:} We proceed by showing
the satisfaction of three conditions in~\cite{Bertsekas1983}: 
(i)~$\gamma$ is sufficiently small, (ii)~$G$ is positive definite, and (iii)~$F$ is positive definite.

\paragraph*{(i)~$\gamma$ is sufficiently small} 
Results in \cite{Bertsekas1991} require
%\begin{equation}
$\gamma \sum_{j=1}^n |H_{ij}| < 1$  for all $i \in \lbrace 1,\ldots, n\rbrace$,
%\end{equation} 
which here follows immediately from~\eqref{eq:gamma}.

\paragraph*{(ii)~$G$ is positive definite}
By definition,~$G$ has only positive diagonal entries. By~$H$'s diagonal dominance we have the following inequality for all~$i \in \lbrace 1,\ldots, n\rbrace$:
\begin{equation}
|G_{ii}| = |H_{ii}| \geq \sum_{\substack{ j=1 \\ j \neq i}}^n |H_{ij}| + \beta > \sum_{\substack{ j=1 \\ j \neq i}}^n |H_{ij}| = \sum_{\substack{ j=1 \\ j \neq i}}^n |G_{ij}|.
\end{equation}
Because~$G$ has positive diagonal entries, is symmetric, and is strictly diagonally dominant,~$G$ is positive definite by Gershgorin's Circle Theorem. 

\paragraph*{(iii)~$F$ is positive definite}
Eq.~\eqref{eq:gamma} ensures the diagonal entries of~$F$ are always positive. And $F$ is diagonally dominant if, for all~$i \in \{1, \ldots, n\}$,
\begin{equation}
|F_{ii}| = 1-\gamma |H_{ii}| > \gamma \sum_{\substack{ j=1 \\ j \neq i}}^n |H_{ij}| = \sum_{\substack{ j=1 \\ j \neq i}}^n |F_{ij}| .
\end{equation}
This requirement can be rewritten as~$\gamma \sum_{j=1}^n |H_{ij}|  < 1$, 
which was satisfied under~(i). 
Because~$F$ has positive diagonal entries, is symmetric, and is strictly diagonally dominant,~$F$ is positive definite by Gershgorin's Circle Theorem.\hfill $\blacksquare$

We next show that the gradient update law~$h$ in~\eqref{eq:grad} converges with asynchronous, distributed computations. Furthermore, we quantify the rate of convergence. 

\begin{lemma} \label{lem:primal}
 \blue{Let~$\gamma$ and~$h$ be as defined in~\eqref{eq:gamma} and~\eqref{eq:grad}.} Let \blue{Assumptions~\ref{as:f}-\ref{as:dualcomm}} hold and fix~$\mu(t) \in \mathcal{M}$. 
Then for the fixed point~$\xhatt$ of~$h$ and for all~$x \in X$, 
\begin{equation}
\Vert h(x)-h(\xhatt) \Vert _\infty \leq q_p\Vert x- \xhatt \Vert _\infty,
\end{equation}
where~$q_p:=(1-\gamma \beta) \in [0,1)$. 
\end{lemma}

\emph{Proof:} For all~$i$, Assumption~\ref{as:diagonal} and the definition of~$F$ give
\begin{equation}
\sum_{\substack{ j=1}}^n F_{ij} = 1- \gamma\Big(|H_{ii}| - \sum_{\substack{ j=1 \\ j \neq i}}^n |H_{ij}|\Big) \leq 1- \gamma\beta.
\end{equation}

This result, \blue{the definition of~$\|\cdot \|_\infty$}, and Lemma~\ref{83lemma} give
\ifbool{Report}{\begin{align*}
\Vert h(x)-h(\xhatt) \Vert _\infty &= \max\limits_i |h_i(x) - h_i(\xhatt)| \\
&\!\leq\!\max\limits_i\!\sum_{\substack{ j=1}}^n\!F_{ij} |x_j\!-\! \hat{x}_{\delta,j}(t)| \\
&\leq\!\max\limits_l\!|x_l\!-\!\hat{x}_{\delta,l}(t)|\!\max\limits_i\!\sum_{\substack{ j=1}}^n\!F_{ij} \\
&\leq \max\limits_l |x_l - \hat{x}_{\delta,l}(t)| (1- \gamma \beta) \\
&= (1-\gamma \beta)\Vert x-\xhatt \Vert _\infty,
\end{align*}}{\begin{align*}
&\Vert h(x)-h(\xhatt) \Vert _\infty = \max\limits_i |h_i(x) - h_i(\xhatt)| \\
&\!\leq\!\max\limits_i\!\sum_{\substack{ j=1}}^n\!F_{ij} |x_j\!-\! \hat{x}_{\delta,j}(t)|\!\leq\!\max\limits_l\!|x_l\!-\!\hat{x}_{\delta,l}(t)|\!\max\limits_i\!\sum_{\substack{ j=1}}^n\!F_{ij} \\
&\leq \max\limits_l |x_l - \hat{x}_{\delta,l}(t)| (1- \gamma \beta) = (1-\gamma \beta)\Vert x-\xhatt \Vert _\infty,
\end{align*}} 
where the last inequality follows from Lemma~\ref{83lemma}. 
All that remains is to show~$(1-\gamma \beta) \in [0,1)$. From~\eqref{eq:gamma} and the inequality~$|H_{ii}|\geq \beta$, for all~$x \in X$ and~$\mu(t) \in M$,
\begin{equation}
\gamma \beta\!<\!\frac{\beta}{\max\limits_i \sum_{j=1}^n |H_{ij}(x, \mu(t))|} \!\leq\!\frac{\beta}{\max\limits_i |H_{ii}(x, \mu(t))|} %\leq \frac{\beta}{\beta} 
\!=\!1. \tag*{$\blacksquare$} 
\end{equation}

%The primal convergence rate then takes the following form.

\begin{lemma} \label{thm:primalconv}
Let \blue{Assumptions~\ref{as:f}-\ref{as:dualcomm}} hold. Let~$\mu(t)$ be the dual vector onboard all primal agents at some time~$k$ and let~$k^t_0$ denote the latest time that any primal agent received the dual variable~$\mu(t)$ that agents currently have onboard. 
Then, with primal agents asynchronously executing
the gradient update law~$h$, agent~$i$ has
\begin{equation}
\|x^i(k;t) - \xhatt \|_\infty \leq q_p^{\textnormal{ops}(k,t)} \max\limits_{j} \|x^j(k^t_0;t) - \xhatt\|_\infty, 
\end{equation}
where~$\xhatt$ is the fixed point of~$h$ with~$\mu(t)$ held constant. 
\end{lemma}
\emph{Proof:} 
From Lemma~\ref{lem:primal} we see that~$h$ is
a~$q_p$-contraction mapping with respect to the norm~$\|\cdot\|_\infty$. 
From Section 6.3 in~\cite{Bertsekas1991}, this
property implies that there
exist sets of the form
\begin{align}
X(k) &= \{x \in \R^n \mid \|x - \xhatt \|_\infty \\
&\leq q_p^k \max_j \|x^j(k^t_0;t ) - \xhatt \|_\infty\}
\end{align}
that satisfy the following criteria from~\cite{hale17}:
\begin{enumerate}[i.]
\item $\cdots \subset X(k+1) \subset X(k) \subset \cdots \subset X$
\item $\lim_{k \to \infty} X(k) = \{\xhatt \}$
\item For all~$i$, there are sets~$X_i(k) \subset X_i$
satisfying
\begin{equation}
X(k) = X_1(k) \times \cdots \times X_N(k)
\end{equation}
\item For all~$y \in X(k)$ and all~$i \in \mathcal{I}_p$, 
$h_i(y) \in X_i(k+1)$, where
%\begin{equation}
$h_i(y) = \Pi_{X_i}\left[y_i - \gamma \nabla_{x_{[i]}} L_{\delta}(y, \mu(t))\right]$. 
%\end{equation}
\end{enumerate}

We will use these properties to compute the desired
convergence rate. Suppose all agents have
a fixed~$\mu(t)$ onboard. Upon receipt of this~$\mu(t)$,
agent~$i$ has~$x^i(k^t_0; t) \in X(0)$ by definition.
Suppose at time~$\ell_i$ that agent~$i$ computes
an update. Then~$x^i_{[i]}(\ell_i+1; t) \in X_i(1)$.
For~$m = \max_{i \in \mathcal{I}_p} \ell_i + 1$, we find
that~$x^i_{[i]}(m; t) \in X_i(1)$ for all~$i$.
Next, suppose that, after all updates have been computed, 
these updated values are sent to and received by all agents
that need them, say at time~$m'$. Then, for any~$i \in \mathcal{I}_p$, agent~$i$
has~$x^i_{[j]}(m'; t) \in X_j(1)$ for all~$j \in \mathcal{I}_p$.
In particular,~$x^i(m'; t) \in X(1)$, and this is satisfied
precisely when ops has incremented by one. Iterating this argument completes the proof. 
\hfill $\blacksquare$

\subsubsection{\blue{Dual Convergence}} \label{sec:dconvergence}
\blue{Towards defining the behavior during an update of a single dual block, we consider the number of operations primal agents compute before communications are sent to a dual agent. In particular, we are interested in defining the oldest primal block a dual agent uses 
in its own computations. 
%Here, ``oldest'' references the primal time~$k$, not necessarily when the update was received by a dual agent. 
%Recall that dual agent~$c$ performs update~$t_c+1$ using its copy of the primal variable,~$x^c(t_c)$. This copy of the primal variable is comprised of updates received from primal agents prior to the time at which dual agent~$c$ computes
%update~$t_c+1$. 
Each of the received primal blocks was sent when some number of operations had been completed by primal agents using the prior dual update. 
Towards quantifying this, we are interested in defining the primal computation time of the oldest primal block used by dual agent~$c$ when
it computes update~$t_c+1$.
\begin{definition} \label{def:kappa}
For dual agent~$c$ computing update~$t_c+1$, let~$\kappa(c,t_c)$ denote the earliest primal computation time for all blocks in~$x^c(t_c)$. That is, for all primal blocks used by dual agent~$c$ during update~$t_c+1$,~$\kappa(c,t_c)$ is the oldest primal time any were computed.% and later sent to dual agent~$c$.
\end{definition}
Thus, the minimum number of operations completed by any primal agent 
for the blocks used by dual agent~$c$ during update~$t_c+1$ is equal to~$\textnormal{ops}(\kappa(c,t_c),t)$.}
We next derive a block-wise
convergence rate for the dual variable.  
%That is, we bound the difference between dual agent~$c$'s updated~$c^{th}$ block, $\mu_{[c]}^c$, and the regularized optimal value for this block, $\mhatc$.

\begin{lemma} \label{thm:dualsingle}
Let Assumptions~\ref{as:f}-\ref{as:dualcomm} hold. 
Let the dual stepsize~$\rho$ be defined such that~$\rho <  \frac{2 \delta}{\delta^2 + 2}$.
Let~$t_c \geq 0$ and consider the case where dual agent~$c$ performs a single update denoted with the iteration
counter~$t_c+1$. Then the distance from the optimal value for block~$c$ is bounded by
\ifbool{Report}{\begin{align*}
\| &\mu^c_{[c]}(t_c+1) - \mhatc \|^2 \leq  q_d \|\mcc - \mhatc \|^2
+ q_p^{2\textnormal{ops}(\kappa(c,t_c),t)} E_1(c) + q_p^{\textnormal{ops}(\kappa(c,t_c),t)} E_2(c) + E_3(c) ,
\end{align*}}{\begin{align*}
\| &\mu^c_{[c]}(t_c+1) - \mhatc \|^2 \leq  q_d \|\mcc - \mhatc \|^2\\
&+ q_p^{2\textnormal{ops}(\kappa(c,t_c),t)} E_1(c) + q_p^{\textnormal{ops}(\kappa(c,t_c),t)} E_2(c) + E_3(c) ,
\end{align*}} 
where~$E_1(c) := (q_d - \rho^2) n M_{[c]}^2 D_{x}^2$,~$E_2(c) := 2 \rho ^2 \sqrt{n} M_{[c]}^2 D_{x}^2$ and~$E_3(c) := (q_d - \rho^2)M_{[c]}^2 D_{x}^2$,~$q_d := (1\!-\!\rho \delta)^2\!+\!2\rho^2 \in [0,1)$,~${M_{[c]} := \max\limits_{x \in X} \| \nabla g_{[c]}(x) \|}$,~$D_{x} := \max\limits_{x,y \in X} \| x-y \|$, 
and~$n$ is the length of the primal variable~$x$.
\end{lemma}

\emph{Proof:} Define~$\xhatt = \argmin_{x \in X} L_{\delta}(x, \mu(t))$ and~$\xhat = \argmin_{x \in X} L_{\delta}(x, \hat{\mu}_\delta)$. Let~$\xctc := x^c(t_c)$ for brevity of notation.
Expanding the dual update law and using the non-expansiveness of~$\Pi_\mathcal{M}$, we find
\ifbool{Report}{\begin{align}
\| \mu^c_{[c]}(t_c\!+\!1) - \mhatc \|^2 &= \|\Pi _{\mathcal{M}_c} [\mcc + \rho (\gc(\xctc)-\delta \mcc)] - \Pi _{\mathcal{M}_c} [\mhatc + \rho (\gc(\xhat)-\delta \mhatc )]\|^2 \\
&\leq  \|\mcc + \rho (\gc(\xctc)-\delta \mcc) - \mhatc  - \rho (\gc(\xhat)-\delta \mhatc)\|^2\ \\
&=  \| (1- \rho \delta)(\mcc - \mhatc) - \rho (\gc(\xhat)-\gc(\xctc))\|^2 \\
&\leq   (1- \rho \delta)^2\|\mcc - \mhatc \|^2_{}  + \rho ^2\|\gc(\xctc) - \gc(\xhat)\|^2_{} \\
&\quad {-2}\rho (1\!-\!\rho \delta)(\mcc\!-\!\mhatc )^T (\gc(\xhat)\!-\!\gc(\xctc)).
\end{align}}{\begin{align}
&\| \mu^c_{[c]}(t_c\!+\!1) - \mhatc \|^2 = \\
&\|\Pi _{\mathcal{M}_c} [\mcc + \rho (\gc(\xctc)-\delta \mcc)]\\
& \quad - \Pi _{\mathcal{M}_c} [\mhatc + \rho (\gc(\xhat)-\delta \mhatc )]\|^2 \\
&\leq  \|\mcc + \rho (\gc(\xctc)-\delta \mcc)\\
& \quad - \mhatc  - \rho (\gc(\xhat)-\delta \mhatc)\|^2\ \\
&=  \| (1- \rho \delta)(\mcc - \mhatc) - \rho (\gc(\xhat)-\gc(\xctc))\|^2 \\
&\leq   (1- \rho \delta)^2\|\mcc - \mhatc \|^2_{}  + \rho ^2\|\gc(\xctc) - \gc(\xhat)\|^2_{} \\
&\quad {-2}\rho (1\!-\!\rho \delta)(\mcc\!-\!\mhatc )^T (\gc(\xhat)\!-\!\gc(\xctc)).
\end{align}} 

Adding~$\gc\left(\xhatt\right) - \gc\left(\xhatt\right)$ inside the last set of parentheses gives
\ifbool{Report}{\begin{align}
\|\mu^c_{[c]}(t_c\!+\!1) - \mhat \|^2 &\leq  (1- \rho \delta)^2\|\mcc - \mhatc \|^2 \nonumber + \rho ^2\|\gc(\xctc) - \gc(\xhat)\|^2 \nonumber \\
& \quad -2 \rho (1-\rho \delta)(\mcc - \mhatc )^T\left(\gc(\xhat)-\gc\left(\xhatt\right)\right) \nonumber \\
& \quad -2\rho (1\!-\!\rho \delta)(\mcc\!-\!\mhatc )^T(\gc\left(\xhatt\right)\!-\!\gc(\xctc)) .
\label{eq:addsub}
\end{align}}{\begin{align}
&\|\mu^c_{[c]}(t_c\!+\!1) - \mhat \|^2 \leq  (1- \rho \delta)^2\|\mcc - \mhatc \|^2 \nonumber \\
& + \rho ^2\|\gc(\xctc) - \gc(\xhat)\|^2 \nonumber \\
& -2 \rho (1-\rho \delta)(\mcc - \mhatc )^T\left(\gc(\xhat)-\gc\left(\xhatt\right)\right) \nonumber \\
& -2\rho (1\!-\!\rho \delta)(\mcc\!-\!\mhatc )^T(\gc\left(\xhatt\right)\!-\!\gc(\xctc)) .
\label{eq:addsub}
\end{align}} 

We can write
\begin{equation*}
0 \leq \|(1-\rho \delta)\left(\gc(\xhat)-\gc\left(\xhatt\right)\right) + \rho (\mcc - \mhatc )\|^2 ,
\end{equation*}
which can be expanded and rearranged to give
\ifbool{Report}{ \begin{align}
&{-2}\rho (1-\rho \delta) (\mcc - \mhatc)^T \left(\gc(\xhat)-\gc\left(\xhatt\right)\right)\\
& \qquad \leq (1\!-\!\rho \delta)^2 \|\gc(\xhat)\!-\!\gc\left(\xhatt\right)\|^2 + \rho ^2 \|\mcc\!-\!\mhatc \|^2 .
\end{align}}{ \begin{align}
&{-2}\rho (1-\rho \delta) (\mcc - \mhatc)^T \left(\gc(\xhat)-\gc\left(\xhatt\right)\right)\\
& \leq (1\!-\!\rho \delta)^2 \|\gc(\xhat)\!-\!\gc\left(\xhatt\right)\|^2 + \rho ^2 \|\mcc\!-\!\mhatc \|^2 .
\end{align}}
Similarly,
\ifbool{Report}{\begin{align}
&{-2} \rho (1-\rho \delta) (\mcc - \mhatc)^T (\gc\left(\xhatt\right)-\gc(\xctc))\\
&\qquad \leq (1\!-\!\rho \delta)^2 \|\gc\left(\xhatt\right)\!-\!\gc(\xctc)\|^2\!+\!\rho ^2 \|\mcc\!-\!\mhatc \|^2 .
\end{align}}{\begin{align}
&{-2} \rho (1-\rho \delta) (\mcc - \mhatc)^T (\gc\left(\xhatt\right)-\gc(\xctc))\\
&\leq (1\!-\!\rho \delta)^2 \|\gc\left(\xhatt\right)\!-\!\gc(\xctc)\|^2\!+\!\rho ^2 \|\mcc\!-\!\mhatc \|^2 .
\end{align}} 
Applying these inequalities to~\eqref{eq:addsub} gives
\ifbool{Report}{\begin{align}
\|\mu^c_{[c]}(t_c+1) - \mhat \|^2 &\leq  (1- \rho \delta)^2\|\mcc - \mhatc \|^2 \\
&\quad +\!\rho ^2\|\gc(\xctc)\!-\!\gc(\xhat)\|^2_{}\!+\!(1\!-\!\rho \delta)^2 \|\gc(\xhat)\!-\!\gc\left(\xhatt\right)\|^2\nonumber \\
& \quad +\!2\rho ^2 \|\mcc\!-\!\mhatc \|^2\!+\!(1\!-\!\rho \delta)^2 \|\gc\left(\xhatt\right)\!-\!\gc(\xctc)\|^2  \nonumber \\
&\leq ((1\!-\!\rho \delta)^2\!+\!2 \rho^2 ) \|\mcc-\mhatc \|^2 +\!\rho ^2\|\gc(\xctc)\!-\!\gc(\xhat)\|^2_{} \nonumber \\
& \quad +(1\!-\!\rho \delta)^2 \|\gc(\xhat)\!-\!\gc\left(\xhatt\right)\|^2 +\!(1\!-\!\rho \delta)^2 \|\gc\left(\xhatt\right)\!-\!\gc(\xctc)\|^2.
\label{eq:gradexp}
\end{align}}{\begin{align}
&\|\mu^c_{[c]}(t_c+1) - \mhat \|^2 \leq  (1- \rho \delta)^2\|\mcc - \mhatc \|^2 \\
&+\!\rho ^2\|\gc(\xctc)\!-\!\gc(\xhat)\|^2_{}\!+\!(1\!-\!\rho \delta)^2 \|\gc(\xhat)\!-\!\gc\left(\xhatt\right)\|^2\nonumber \\
&+\!2\rho ^2 \|\mcc\!-\!\mhatc \|^2\!+\!(1\!-\!\rho \delta)^2 \|\gc\left(\xhatt\right)\!-\!\gc(\xctc)\|^2  \nonumber \\
&\leq ((1\!-\!\rho \delta)^2\!+\!2 \rho^2 ) \|\mcc-\mhatc \|^2 \nonumber \\
&+\!\rho ^2\|\gc(\xctc)\!-\!\gc(\xhat)\|^2_{}\!+\!(1\!-\!\rho \delta)^2 \|\gc(\xhat)\!-\!\gc\left(\xhatt\right)\|^2 \nonumber \\
&+\!(1\!-\!\rho \delta)^2 \|\gc\left(\xhatt\right)\!-\!\gc(\xctc)\|^2.
\label{eq:gradexp}
\end{align}} 

%Looking to manipulate the~$\rho^2\|\gc(\xctc) - \gc(\xhat)\|^2$ term, we can write
%\begin{align}
%&\|\gc(\xctc)\!-\!\gc(\xhat)\|^2\!=\! \|\gc(\xctc)\!-\!\gc(\xhatt)\!+\!\gc(\xhatt)\!-\!\gc(\xhat)\|^2 \nonumber \\
%&\leq \|\gc(\xctc) - \gc(\xhatt)\|^2 + \|\gc(\xhatt) - \gc(\xhat)\|^2 \nonumber \\
%&+ 2\|\gc(\xctc) - \gc(\xhatt)\| \|\gc(\xhatt) - \gc(\xhat)\|,
%\end{align}
%and substituting into~\eqref{eq:gradexp} gives
In~\eqref{eq:gradexp}, we next use~$\rho^2\|g_{[c]}(x^c_t) - g_{[c]}(\hat{x}_{\delta})\|^2 = \rho^2\|g_{[c]}(x^c_t) - g_{[c]}\left(\xhatt\right) + g_{[c]}\left(\xhatt\right) - g_{[c]}(\hat{x}_{\delta})\|^2$, then expand, 
and combine like terms to find
%in gives
%\begin{align}
%&\| \mu^c_{[c]}(t_c+1) - \mhatc \|^2 \leq  ((1\!-\!\rho \delta)^2\!+\!2\rho^2) \|\mcc\!-\!\mhatc \|^2 \nonumber \\
%& + \rho ^2\|\gc(\xctc) - \gc\left(\xhatt\right)\|^2 + \rho^2 \|\gc\left(\xhatt\right) - \gc(\xhat)\|^2 \nonumber \\
%& + 2 \rho ^2 \|\gc(\xctc) - \gc\left(\xhatt\right)\| \|\gc\left(\xhatt\right) - \gc(\xhat)\| \nonumber \\
%& + (1-\rho \delta)^2 \|\gc(\xhat)-\gc\left(\xhatt\right)\|^2 \nonumber \\
%& + (1-\rho \delta)^2 \|\gc\left(\xhatt\right)-\gc(\xctc)\|^2  .
%\end{align}
%Grouping terms, this inequality simplifies to
\ifbool{Report}{\begin{align}
\| \mu^c_{[c]}(t_c+1) - \mhatc \|^2 &\leq  ((1\!-\!\rho \delta)^2\!+\!2\rho^2) \|\mcc\!-\!\mhatc \|^2 + ((1-\rho\delta)^2 + \rho ^2)\|\gc(\xctc) - \gc\left(\xhatt\right)\|^2  \nonumber \\
&\quad + 2 \rho ^2 \|\gc(\xctc) - \gc\left(\xhatt\right)\| \|\gc\left(\xhatt\right) - \gc(\xhat)\| \nonumber \\
&\quad + ((1-\rho\delta)^2 + \rho ^2) \|\gc\left(\xhatt\right) - \gc(\xhat)\|^2  .
\end{align}}{\begin{align}
&\| \mu^c_{[c]}(t_c+1) - \mhatc \|^2 \leq  ((1\!-\!\rho \delta)^2\!+\!2\rho^2) \|\mcc\!-\!\mhatc \|^2 \nonumber \\
&\quad + ((1-\rho\delta)^2 + \rho ^2)\|\gc(\xctc) - \gc\left(\xhatt\right)\|^2  \nonumber \\
&\quad + 2 \rho ^2 \|\gc(\xctc) - \gc\left(\xhatt\right)\| \|\gc\left(\xhatt\right) - \gc(\xhat)\| \nonumber \\
&\quad + ((1-\rho\delta)^2 + \rho ^2) \|\gc\left(\xhatt\right) - \gc(\xhat)\|^2  .
\end{align}} 
Using the Lipschitz property of~$g_{[c]}$, we can write
\ifbool{Report}{\begin{align}
\| \mu^c_{[c]}(t_c+1) - \mhatc \|^2 &\leq  ((1\!-\!\rho \delta)^2\!+\!2\rho^2) \|\mcc\!-\!\mhatc \|^2 + ((1-\rho\delta)^2 + \rho ^2)M_{[c]}^2\|\xctc - \xhatt\|^2  \nonumber \\
&\quad + 2 \rho ^2 M_{[c]}^2 \|\xctc - \xhatt\| \|\xhatt - \xhat\| + ((1-\rho\delta)^2 + \rho ^2)M_{[c]}^2 \|\xhatt - \xhat\|^2  .
\end{align}}{\begin{align}
&\| \mu^c_{[c]}(t_c+1) - \mhatc \|^2 \leq  ((1\!-\!\rho \delta)^2\!+\!2\rho^2) \|\mcc\!-\!\mhatc \|^2 \nonumber \\
&\quad + ((1-\rho\delta)^2 + \rho ^2)M_{[c]}^2\|\xctc - \xhatt\|^2  \nonumber \\
&\quad + 2 \rho ^2 M_{[c]}^2 \|\xctc - \xhatt\| \|\xhatt - \xhat\| \nonumber \\
&\quad + ((1-\rho\delta)^2 + \rho ^2)M_{[c]}^2 \|\xhatt - \xhat\|^2  .
\end{align}} 
Using~$\|\xhatt - \xhat\| \leq D_{x}$, the inequality simplifies to
\ifbool{Report}{\begin{align}
\| \mu^c_{[c]}(t_c+1) - \mhatc \|^2 &\leq  ((1\!-\!\rho \delta)^2\!+\!2\rho^2) \|\mcc\!-\!\mhatc \|^2 + ((1\!-\!\rho\delta)^2\!+\!\rho ^2)M_{[c]}^2\|\xctc\!-\!\xhatt\|^2\nonumber \\
& \quad + 2 \rho ^2 M_{[c]}^2 D_{x} \|\xctc\!-\!\xhatt\|  + ((1\!-\!\rho\delta)^2\!+\!\rho ^2)M_{[c]}^2 D_{x}^2  .
\label{eq:Dmax}
\end{align}}{\begin{align}
&\| \mu^c_{[c]}(t_c+1) - \mhatc \|^2 \leq  ((1\!-\!\rho \delta)^2\!+\!2\rho^2) \|\mcc\!-\!\mhatc \|^2 \nonumber \\
& + ((1\!-\!\rho\delta)^2\!+\!\rho ^2)M_{[c]}^2\|\xctc\!-\!\xhatt\|^2 + 2 \rho ^2 M_{[c]}^2 D_{x} \|\xctc\!-\!\xhatt\| \nonumber \\
& + ((1\!-\!\rho\delta)^2\!+\!\rho ^2)M_{[c]}^2 D_{x}^2  .
\label{eq:Dmax}
\end{align}} 

Using Definition~\ref{def:kappa}, define~$\tilde{x}^c(t_c)$ as the primal variable whose distance is greatest from the optimal value at primal time~$\kappa(c,t_c)$. 
That is,~$\tilde{x}^c(t_c) := \max_{j \in \mathcal{I}_p} \|x^j(\kappa(c,t_c), t_c) - \xhatt \|$. 
Using this, the contraction property of primal updates from \blue{Lemma}~\ref{thm:primalconv}, and the definition of~$D_x$, we find
\ifbool{Report}{\begin{align}
    \|\xctc - \xhatt \| &\leq \| \tilde{x}^c(t_c) - \xhatt \| \leq \sqrt{n} \| \tilde{x}^c(t_c) - \xhatt \|_\infty 
    %&\leq q_p^{\textnormal{ops}(\kappa(c,t_c),t)} \sqrt{n} \max_j \|x^j(k^t_0; t)- \xhatt \|_\infty \\
    \leq q_p^{\textnormal{ops}(\kappa(c,t_c),t)} \sqrt{n} D_x.
    \label{eq:primalnorm}
\end{align}}{\begin{align}
    \|\xctc - \xhatt \| &\leq \| \tilde{x}^c(t_c) - \xhatt \| \leq \sqrt{n} \| \tilde{x}^c(t_c) - \xhatt \|_\infty \\
    %&\leq q_p^{\textnormal{ops}(\kappa(c,t_c),t)} \sqrt{n} \max_j \|x^j(k^t_0; t)- \xhatt \|_\infty \\
    &\leq q_p^{\textnormal{ops}(\kappa(c,t_c),t)} \sqrt{n} D_x.
    \label{eq:primalnorm}
\end{align}} 

Applying this result to~\eqref{eq:Dmax} above gives
\ifbool{Report}{\begin{align}
\| \mu^c_{[c]}(t_c+1) - \mhatc \|^2 &\leq  ((1\!-\!\rho \delta)^2\!+\!2\rho^2) \|\mcc\!-\!\mhatc \|^2 + ((1\!-\!\rho\delta)^2\!+\!\rho ^2) n M_{[c]}^2 q_p^{2\textnormal{ops}(\kappa(c,t_c),t)} D_{x}^2  \nonumber \\
&\quad + 2 \rho ^2 \sqrt{n} M_{[c]}^2 D_{x}^2 q_p^{\textnormal{ops}(\kappa(c,t_c),t)}  + ((1\!-\!\rho\delta)^2\!+\!\rho ^2)M_{[c]}^2 D_{x}^2 .
\end{align}}{\begin{align}
&\| \mu^c_{[c]}(t_c+1) - \mhatc \|^2 \leq  ((1\!-\!\rho \delta)^2\!+\!2\rho^2) \|\mcc\!-\!\mhatc \|^2 \nonumber \\
&\quad + ((1\!-\!\rho\delta)^2\!+\!\rho ^2) n M_{[c]}^2 q_p^{2\textnormal{ops}(\kappa(c,t_c),t)} D_{x}^2  \nonumber \\
&\quad + 2 \rho ^2 \sqrt{n} M_{[c]}^2 D_{x}^2 q_p^{\textnormal{ops}(\kappa(c,t_c),t)}  + ((1\!-\!\rho\delta)^2\!+\!\rho ^2)M_{[c]}^2 D_{x}^2 .
\end{align}} 
Using~$\rho< \frac{2 \delta}{\delta^2 + 2}$, we have~$q_d = (1-\rho\delta)^2 + 2 \rho^2 \in (0,1)$, completing the proof.
\hfill $\blacksquare$

\begin{lemma} \label{lem:mufinal}
\blue{Let all conditions and definitions of \blue{Lemma}~\ref{thm:dualsingle} hold.
Let~$T(t) = \min_c t_c$ be the minimum number of updates any one dual agent has performed by time~$t$ and let~$K(t)$ be the minimum number of operations primal agents completed on any primal block used to compute any dual block from~$\mu(0)$ to~$\mu(t)$. Then, Algorithm~\ref{alg:2}'s convergence for~$\mu$ obeys
\ifbool{Report}{\begin{align*}
    \| \mu(t)\!-\!\mhatd \|^2\!&\leq\!q_d^{T(t)} \|\mu (0)\!-\!\mhatd \|^2\!+\!\Bigl(\!q_p^{2K(t)} (q_d\!-\!\rho^2) n N_d M^2 D_{x}^2 \\
    &\qquad + q_p^{K(t)} 2 \rho ^2 \sqrt{n} N_d M^2 D_{x}^2+ (q_d - \rho^2) N_d M^2 D_{x}^2 \Bigr) \frac{1}{1-q_d} ,
\end{align*}}{\begin{align*}
    &\| \mu(t)\!-\!\mhatd \|^2\!\leq\!q_d^{T(t)} \|\mu (0)\!-\!\mhatd \|^2\!+\!\Bigl(\!q_p^{2K(t)} (q_d\!-\!\rho^2) n N_d M^2 D_{x}^2 \\
    &\quad + q_p^{K(t)} 2 \rho ^2 \sqrt{n} N_d M^2 D_{x}^2+ (q_d - \rho^2) N_d M^2 D_{x}^2 \Bigr) \frac{1}{1-q_d} ,
\end{align*}} 
where~$M := \max\limits_{x \in X} \| \nabla g(x) \|$ and~$N_d$ is the number of dual agents.}
\end{lemma}
\blue{\emph{Proof:}
Let~$K_c(t_c)$ be the minimum number of operations primal agents completed on any primal block used to compute~$\mu_{[c]}$ from~$\mu_{[c]}(0)$ to~$\mu_{[c]}(t_c)$. Then recursively applying \blue{Lemma}~\ref{thm:dualsingle} and using the definition of~$K_c(t_c)$ gives
\ifbool{Report}{\begin{align}
    \| \mu_{[c]}^c (t_c) \!-\! \mhat \|^2 \!&\leq\! q_d \|\mu_{[c]}^c (t_c-1) - \mhatc \|^2 + q_p^{2K_c(t_c)} E_1(c) + q_p^{K_c(t_c)} E_2(c) + E_3(c) \\
%& \leq q_d^2 \|\mu_{[c]}^c (a-1) - \mhatc \|^2 + q_d q_p^{2\textnormal{ops}(\kappa(c,a\!-\!1),T^c_{a\!-\!1})} E_1(T^c_{a\!-\!1}) \\
%&\quad + q_d q_p^{\textnormal{ops}(\kappa(c,a\!-\!1),T^c_{a\!-\!1})} E_2(T^c_{a\!-\!1})+ q_d E_3 \\
%&\quad + q_p^{2\textnormal{ops}(\kappa(c,a),T^c_a)} E_1(T^c_a) + q_p^{\textnormal{ops}(\kappa(c,a),T^c_a)} E_2(T^c_a) + E_3\\
&\!\leq q_d^{t_c} \|\mu_{[c]}^c (0)\!-\!\mhatc \|^2\!+\!\sum_{i=0}^{t_c - 1}\!q_d^i \Big( q_p^{2K_c(t_c)} E_1(c)  +  q_p^{K_c(t_c)} E_2(c) +  E_3(c) \Big) \\
&\leq q_d^{t_c} \|\mu_{[c]}^c (0) - \mhatc \|^2 +  \Bigl( q_p^{2K_c(t_c)} E_1(c) +  q_p^{K_c(t_c)} E_2(c) +  E_3(c)\Bigr)\frac{1-q_d^{t_c}}{1-q_d}, \label{eq:dualgeometric}
    \end{align}}{\begin{align}
    &\| \mu_{[c]}^c (t_c) \!-\! \mhat \|^2 \!\leq\! q_d \|\mu_{[c]}^c (t_c-1) - \mhatc \|^2 \\
    &\quad + q_p^{2K_c(t_c)} E_1(c) + q_p^{K_c(t_c)} E_2(c) + E_3(c) \\
%& \leq q_d^2 \|\mu_{[c]}^c (a-1) - \mhatc \|^2 + q_d q_p^{2\textnormal{ops}(\kappa(c,a\!-\!1),T^c_{a\!-\!1})} E_1(T^c_{a\!-\!1}) \\
%&\quad + q_d q_p^{\textnormal{ops}(\kappa(c,a\!-\!1),T^c_{a\!-\!1})} E_2(T^c_{a\!-\!1})+ q_d E_3 \\
%&\quad + q_p^{2\textnormal{ops}(\kappa(c,a),T^c_a)} E_1(T^c_a) + q_p^{\textnormal{ops}(\kappa(c,a),T^c_a)} E_2(T^c_a) + E_3\\
&\!\leq q_d^{t_c} \|\mu_{[c]}^c (0)\!-\!\mhatc \|^2\!+\!\sum_{i=0}^{t_c - 1}\!q_d^i \Big( q_p^{2K_c(t_c)} E_1(c) \\
&\quad +  q_p^{K_c(t_c)} E_2(c) +  E_3(c) \Big) \\
&\leq q_d^{t_c} \|\mu_{[c]}^c (0) - \mhatc \|^2 +  \Bigl( q_p^{2K_c(t_c)} E_1(c) \\
    &\quad +  q_p^{K_c(t_c)} E_2(c) +  E_3(c)\Bigr)\frac{1-q_d^{t_c}}{1-q_d}, \label{eq:dualgeometric}
    \end{align}} 
where the last inequality uses~$q_d \in [0,1)$ and sums the geometric series. We now derive a bound on the entire~$\mu$ vector at time~$t$. Expanding~$\| \mu(t) - \mhatd \|^2$ allows us to write
\ifbool{Report}{\begin{align}
    \| \mu(t) - \mhatd \|^2 &= \sum_{c=1}^{N_d} \| \mu_{[c]}^c (t_c) \!-\! \mhat \|^2 \\
    &\leq \sum_{c=1}^{N_d} q_d^{t_c} \|\mu_{[c]}^c (0) - \mhatc \|^2 +  \Bigl( q_p^{2K_c(t_c)} E_1(c) +  q_p^{K_c(t_c)} E_2(c) +  E_3(c)\Bigr)\frac{1-q_d^{t_c}}{1-q_d} \\
    &\leq \sum_{c=1}^{N_d} q_d^{t_c} \|\mu_{[c]}^c (0)\!-\!\mhatc \|^2 +  \Bigl( q_p^{2K(t)}(q_d\!-\!\rho^2) n N_d M^2 D_{x}^2 \\
    &\qquad +  q_p^{K(t)} 2 \rho ^2 \sqrt{n} N_d M^2 D_{x}^2 +  (q_d - \rho^2)N_d M^2 D_{x}^2\Bigr)\frac{1}{1-q_d},
\end{align}}{\begin{align}
    &\| \mu(t) - \mhatd \|^2 = \sum_{c=1}^{N_d} \| \mu_{[c]}^c (t_c) \!-\! \mhat \|^2 \\
    &\leq \sum_{c=1}^{N_d} q_d^{t_c} \|\mu_{[c]}^c (0) - \mhatc \|^2 +  \Bigl( q_p^{2K_c(t_c)} E_1(c) \\
    &\quad +  q_p^{K_c(t_c)} E_2(c) +  E_3(c)\Bigr)\frac{1-q_d^{t_c}}{1-q_d} \\
    &\leq \sum_{c=1}^{N_d} q_d^{t_c} \|\mu_{[c]}^c (0)\!-\!\mhatc \|^2 +  \Bigl( q_p^{2K(t)}(q_d\!-\!\rho^2) n N_d M^2 D_{x}^2 \\
    &\quad +  q_p^{K(t)} 2 \rho ^2 \sqrt{n} N_d M^2 D_{x}^2 +  (q_d - \rho^2)N_d M^2 D_{x}^2\Bigr)\frac{1}{1-q_d},
\end{align}} }
\blue{where the first inequality applies~\eqref{eq:dualgeometric} and the second uses~$K_c(t_c) \geq K(t)$,~$M_{[c]}^2 \leq N_d M^2$, and simplifies. Applying the summation and definition of~$T(t)$ completes the proof. \hfill $\blacksquare$}

\subsubsection{\blue{Proof of Theorem~\ref{thm:final}}} 
\blue{We see that 
\ifbool{Report}{\begin{align}
\|x^i(k; t) \!-\!   \hat{x}_{\delta}\|^2 &= \|x^i(k; t) \!-\!\xhatt + \xhatt \!-\!  \hat{x}_{\delta}\|^2 \\
&\leq 2\|x^i(k; t) -\xhatt \|^2 +  2\|\xhatt - \hat{x}_{\delta}\|^2 \\
&\leq 2 n \|x^i(k; t) -\xhatt \|^2_\infty +  \frac{2 M^2} {\beta^2} \|\mu(t) - \hat{\mu}_{\delta}\|^2,
\end{align}}{\begin{align}
&\|x^i(k; t) \!-\!   \hat{x}_{\delta}\|^2 = \|x^i(k; t) \!-\!\xhatt + \xhatt \!-\!  \hat{x}_{\delta}\|^2 \\
&\leq 2\|x^i(k; t) -\xhatt \|^2 +  2\|\xhatt - \hat{x}_{\delta}\|^2 \\
&\leq 2 n \|x^i(k; t) -\xhatt \|^2_\infty +  \frac{2 M^2} {\beta^2} \|\mu(t) - \hat{\mu}_{\delta}\|^2,
\end{align}} 
where the last line applies Lemma~4.1 in~\cite{koshal2011multiuser}. Next, 
applying Lemmas~\ref{thm:primalconv} and~\ref{lem:mufinal} gives
\ifbool{Report}{\begin{align}
    \|x^i(k; t)\!- \! \hat{x}_{\delta}\|^2\!&\leq\!2nq_p^{2 \textnormal{ops}(k,t)}\!\max\limits_{j} \!\|x^j\!(k^t_0;t)\!-\!\xhatt\|^2_\infty \\
    &\qquad + q_d^{T(t)} \frac{2 M^2} {\beta^2} \|\mu (0)\!-\!\mhatd \|^2  + \left( q_p^{2K(t)}  \frac{2 n N_d M^4 D_{x}^2 (q_d\!-\!\rho^2)} {\beta^2} \right. \\
    &\qquad + q_p^{K(t)} \frac{4 \rho^2 \sqrt{n} N_d M^4 D_{x}^2} {\beta^2} +  \left. \frac{2 N_d  M^4 D_{x}^2 (q_d - \rho^2)} {\beta^2} \right) \frac{1}{1-q_d}.
\end{align}}{\begin{align}
    &\|x^i(k; t)\!- \! \hat{x}_{\delta}\|^2\!\leq\!2nq_p^{2 \textnormal{ops}(k,t)}\!\max\limits_{j} \!\|x^j\!(k^t_0;t)\!-\!\xhatt\|^2_\infty \\
    &+ q_d^{T(t)} \frac{2 M^2} {\beta^2} \|\mu (0)\!-\!\mhatd \|^2  + \left( q_p^{2K(t)}  \frac{2 n N_d M^4 D_{x}^2 (q_d\!-\!\rho^2)} {\beta^2} \right. \\
    &+ q_p^{K(t)} \frac{4 \rho^2 \sqrt{n} N_d M^4 D_{x}^2} {\beta^2} +  \left. \frac{2 N_d  M^4 D_{x}^2 (q_d - \rho^2)} {\beta^2} \right) \frac{1}{1-q_d}.
\end{align}} 
Defining~$C_1,$~$C_2,$ and~$C_3$ completes the proof. \hfill $\blacksquare$}

\subsection{\blue{Proof of Corollary~\ref{cor:asynchpen}:}} \label{app:asynchpen}
\blue{We first simplify by noting that~$2 \textnormal{ops}(k,t) \geq K(t)$ and~$2K(t) \geq K(t)$. This allows us to 
factor the bound in Theorem~\ref{thm:final} with~$q_p^{K(t)} \Bigl( 2 n \max\limits_{j} \|x^j(k^t_0;t)\!-\!\xhatt\|^2_\infty +  C_1 + C_2 \Bigr).$
% \begin{align}
%     q_p^{K(t)} \Bigl( 2 n \max\limits_{j} \|x^j(k^t_0;t)\!-\!\xhatt\|^2_\infty +  C_1 + C_2 \Bigr).
% \end{align}
Setting this less than or equal to~$\frac{\epsilon_1}{2}$ and solving gives the lower bound on~$K(t)$.
% \begin{align}
%     K(t)\!\geq\!\frac{\log(\epsilon_1)\!-\!\log(4n D_x^2 \!+\!2C_1\!+\!2C_2)}{\log(q_p)}.
% \end{align}
Similarly, setting~$q_d^{T(t)} \frac{2 M^2} {\beta^2} \|\mu (0) - \mhatd \|^2 \leq \frac{\epsilon_1}{2}$ gives the lower bound on~$T(t)$.
% \begin{align}
%     T(t)  \geq \frac{\log (\epsilon_1 \beta^2) - \log(4 M^2 \|\mu (0) - \mhatd \|^2)}{\log (q_d)}.
% \end{align}
Finally, we set~$\rho = \frac{\delta}{1+\delta^2}$ which results in~$(q_d - \rho^2) = \frac{1}{1+\delta^2}$.
% \begin{align}
%     (q_d - \rho^2) &= (1-\rho \delta)^2 + \rho^2 \\
%     &= (1 - \frac{\delta^2}{1+\delta^2})^2 + \frac{\delta^2}{(1+\delta^2)^2} \\
%     &= \frac{1+ \delta^2 }{(1+\delta^2)^2} \\
%     &= \frac{1}{1+\delta^2}.
% \end{align}
Applying this to~$C_3$ and setting less than or equal to~$\epsilon_2$ gives the final bound on~$\delta^2$.
% \begin{align}
%     \frac{1}{1+\delta^2} \frac{2 N_d  M^4 D_{x}^2} {\beta^2 (1-q_d)} &\leq \epsilon_2 \\
%     \frac{1}{1+\delta^2} &\leq \frac{\epsilon_2 \beta^2 (1-q_d) }{2 N_d  M^4 D_{x}^2} \\
%     \delta^2 &\geq \frac{2 N_d  M^4 D_{x}^2}{\epsilon_2 \beta^2 (1-q_d) } - 1.
% \end{align}
 \hfill $\blacksquare$}

\bibliographystyle{plain}
\bibliography{sources2}

\ifbool{Report}{}{\vspace{-0.8cm} 

\begin{IEEEbiography}[{\includegraphics[width=1in,height=1.25in,clip,keepaspectratio,draft=false]{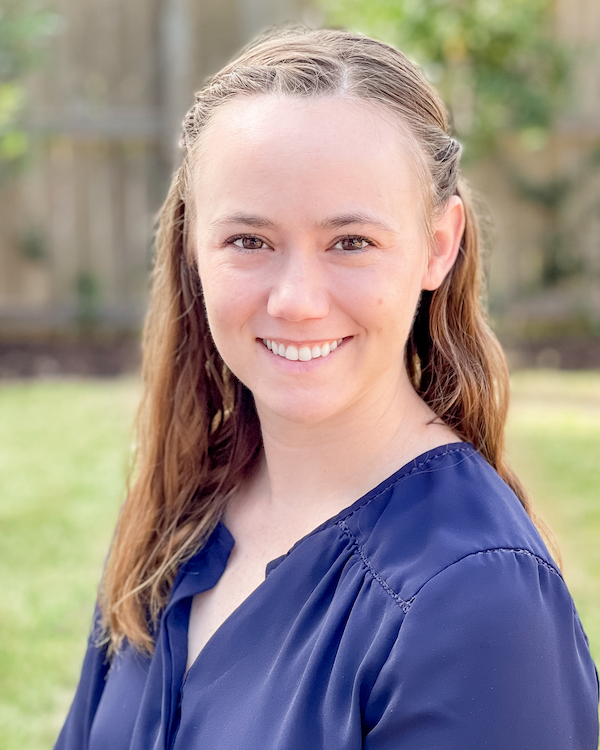}}] {Katherine Hendrickson} is a PhD candidate at the University of Florida. She received her BS in Mathematics from Auburn University in 2013 and MS in Industrial and Systems Engineering from the University of Florida in 2017. She has worked professionally in the field of survivability and lethality. Her research interests include constrained optimization and hybrid systems with an emphasis on asynchronous multi-agent systems and distributed optimization.
\end{IEEEbiography}

\vspace{-0.8cm} 

\begin{IEEEbiography}[{\includegraphics[width=1in,height=1.25in,clip,keepaspectratio,draft=false]{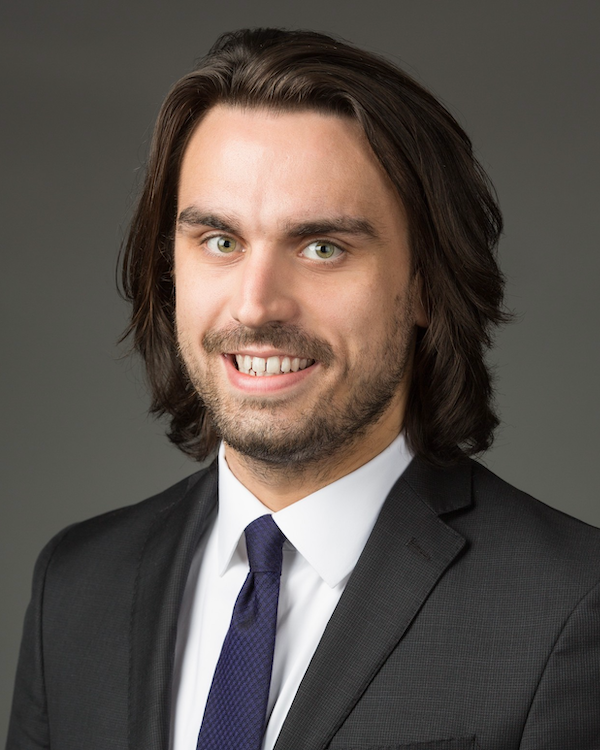}}] {Matthew Hale} is an Assistant Professor of Mechanical and Aerospace Engineering at the University of Florida. He received his BSE in Electrical Engineering \emph{summa cum laude} from the University of Pennsylvania in 2012, and his MS and PhD in Electrical and Computer Engineering from the Georgia Institute of Technology in 2015 and 2017, respectively. His research interests include multi-agent systems, mobile robotics, privacy in control, and distributed optimization. 
He received an NSF CAREER Award in 2020 and an ONR YIP in 2022.   
\end{IEEEbiography}}

\end{document}